\documentclass[12pt]{amsart}

\usepackage[T1]{fontenc}
\usepackage[utf8]{inputenc}
\usepackage{a4wide}
\usepackage{amsmath,amssymb,amsthm,mathtools}
\usepackage{xcolor}
\usepackage{titlesec}
\usepackage{hyperref}

\titleformat{\section}
 {\color{blue}\normalfont\Large\bfseries}
 {\thesection}{1em}{}
\hypersetup{colorlinks=true,linkcolor=blue,citecolor=cyan,urlcolor=blue}
\makeatletter
\renewcommand{\@pnumwidth}{2em}
\def\l@section{\@tocline{1}{0pt}{1pc}{3pc}{}}
\makeatother
\newtheorem{theorem}{Theorem}
\newtheorem{proposition}{Proposition}
\newtheorem{lemma}{Lemma}
\newtheorem{corollary}{Corollary}
\theoremstyle{remark}
\newtheorem{remark}{Remark}

\newcommand{\F}{\mathcal F}
\newcommand{\Fplus}{\mathcal F_{+}^{(1)}}
\newcommand{\M}{\mathcal M}
\newcommand{\ind}{\mathbf 1}
\newcommand{\eps}{\varepsilon}
\newcommand{\dd}{\,\mathrm d}
\newcommand{\e}{\mathrm e}
\newcommand{\Ssmooth}{\mathcal S}

\allowdisplaybreaks[3]

\begin{document}

\title[Large values of derivatives of quadratic Dirichlet $L$-functions]
 {Large values of derivatives of quadratic Dirichlet $L$-functions}
\author{Zikang Dong}
\author{Ruihua Wang}
\author{Weijia Wang}
\author{Hao Zhang}
\author{Shengbo Zhao}
\address[Zikang Dong]{School of Mathematical Sciences, Soochow University,
 Suzhou 215006, P. R. China}
\address[Ruihua Wang]{School of Fundamental Sciences, Hainan Bielefeld
 University of Applied Sciences, Danzhou 578101, P. R. China}
\address[Weijia Wang]{Morningside Center of Mathematics, Academy of
 Mathematics and Systems Science, Chinese Academy of Sciences, Beijing
 100190, P. R. China}
\address[Hao Zhang]{School of Mathematics, Hunan University, Changsha
 410082, P. R. China}
\address[Shengbo Zhao]{School of Mathematical Sciences, Key Laboratory of
 Intelligent Computing and Applications (Tongji University), Ministry of
 Education, Tongji University, Shanghai 200092, China}
\email{zikangdong@gmail.com}
\email{ruih.wan9@gmail.com}
\email{weijiawang@amss.ac.cn}
\email{zhanghaomath@hnu.edu.cn}
\email{shengbozhao@hotmail.com}
\date{October 13, 2026}

\begin{abstract}
We study large values of fixed derivatives of quadratic Dirichlet
$L$-functions on the real segment $1/2\leq\sigma\leq1$.  Assuming GRH,
we prove lower bounds in terms of GCD sums and finite Euler products.
At the central point, every fixed derivative order is compared with the
same square-free GCD sum.  In the mesoscopic band
$\sigma=1/2+A/\log_2X$, a renormalized resonator estimates contributions
from both halves of the differentiated approximate functional equation.
The resulting lower bound has a relative error.  On a smaller boundary
scale, a common finite product appears in both contributions to this
lower bound.
This allows us to locate the crossover between their derivative weights.
The coarser exponential bounds have constant $1$ at the center and
$\mathrm e^{-A}$ in the mesoscopic band.  For fixed $1/2<\sigma<1$,
the lower bound includes the full typical logarithmic weight of the
resonating Euler product.  Burgess truncation improves the exponential
constant when a shorter cutoff allows a longer resonator.
Near the $1$-line, our lower bounds involve Dickman constants
on the $1/\log_2X$ scale.  We also prove a uniform transition estimate
outside a shrinking neighborhood of the central point.
All results hold in the positive even subfamily of fundamental
discriminants.
\end{abstract}

\maketitle
\tableofcontents

\section{Introduction}

For $j\geq1$, $\log_j$ denotes the $j$-fold iterated natural logarithm.
Throughout, the derivative order $\ell$ is a fixed nonnegative integer.
Let $\F$ be the set of fundamental discriminants, and write
$\chi_d=(\frac d\cdot)$ for the primitive quadratic character of conductor
$|d|$.  We shall work in the subfamily
\[
 \Fplus:=\{d>0:d\text{ is square-free and }d\equiv1\pmod4\}.
\]
Every member of $\Fplus$ is a positive fundamental discriminant, and its
character is even.  Lower bounds for the maximum over $\Fplus$
also hold over $\F$.
For the local densities used below, put
\[
 h(n):=\prod_{\substack{p\mid n\\p>2}}\frac p{p+1}.
\]
We write $P^+(n)$ for the largest prime factor of $n$, with $P^+(1)=1$.

\medskip
\noindent\textbf{Character sums and their role.}
For a nonprincipal Dirichlet character $\chi\pmod q$, write
\[
 S_\chi(M,H):=\sum_{M<n\leq M+H}\chi(n),\qquad 1\leq H\leq q.
\]
The P\'olya--Vinogradov inequality gives, uniformly in $M$ and $H$,
\begin{equation}\label{eq:PV-background}
 |S_\chi(M,H)|\ll\sqrt q\log q.
\end{equation}
When $q$ is cube-free, Burgess's bound is available for every fixed
integer $r\geq1$:
\begin{equation}\label{eq:Burgess-background}
 |S_\chi(M,H)|\ll_{r,\varepsilon}
 H^{1-1/r}q^{(r+1)/(4r^2)+\varepsilon}.
\end{equation}
The conductors in $\Fplus$ are square-free, so all these choices of
$r$ are available here.  Standard accounts appear in
\cite[Chapter~9]{MV} and \cite[Chapter~12]{IK}.
Refinements of P\'olya--Vinogradov relate large character sums to
small-conductor characters \cite{GSchar,GSPV}.
Under GRH, the approximation of long character sums by sums over
integers with small prime factors also leads to the Dickman function
near the $1$-line \cite{GSchar}.

These estimates serve two different purposes in this paper.
First, partial summation of \eqref{eq:PV-background} or
\eqref{eq:Burgess-background} truncates the series
\[
 (-1)^\ell L^{(\ell)}(\sigma,\chi_d)
 =\sum_{n\geq1}\frac{\chi_d(n)(\log n)^\ell}{n^\sigma}.
\]
A shorter Burgess cutoff permits a longer resonator and leads to the
optimization over $r$ in Corollary~\ref{cor:fixed-exponent}.
Second, expanding a resonator moment produces a second type of sum,
$\sum_d\chi_d(n)$, over discriminants.
Quadratic reciprocity and GRH give a positive main term for square
$n$ and a controlled error otherwise. See Lemma~\ref{lem:orth}.
Burgess's estimate is unconditional here.
GRH enters through this orthogonality estimate for the family.

Soundararajan introduced the modern resonance method \cite{SoundRes}, and
long resonators later connected central values with GCD sums
\cite{BondarenkoSeip,BT}.  The fixed-strip problem for all characters
modulo a prime, and recently for their derivatives, is treated in
\cite{AMMP,CuiPengSongZhao}.  Exact character-group orthogonality is not
available in the quadratic subfamily.  Controlling the nonsquare terms
in the orthogonality relation is the additional difficulty addressed
by Darbar and Maiti \cite{DarbarMaiti}.

Our main theorems compare derivative weights with specified GCD sums
or finite Euler products.  We state the leading exponential consequences
separately. See Remark~\ref{rem:weight-precision} for this distinction.

The central-point transfer holds for any square-free set with the
following prime-support bound. A sharp GCD-sum construction will then
give its exponential consequence.

\begin{theorem}\label{thm1}
Assume GRH, fix an integer $\ell\geq0$ and $0<\beta<1/4$, and put
$N=\lfloor X^{1/4-\beta}\rfloor$.  Let $\mathcal M$ be a set of $N$
distinct square-free positive integers satisfying
\[
 \max_{m\in\mathcal M}P^+(m)\leq(\log N)^{1+o(1)}.
\]
There is a constant $c_\ell>0$, depending only on $\ell$, such that
\begin{equation}\label{eq:central-precise}
 \max_{\substack{X<d\leq2X\\d\in\Fplus}}
 (-1)^\ell L^{(\ell)}(\tfrac12,\chi_d)
 \geq(c_\ell+o(1))(\log X)^\ell
 \frac1N\sum_{m,n\in\mathcal M}
 \sqrt{\frac{(m,n)}{[m,n]}}.
\end{equation}
The error tends to zero for each family of sets satisfying the displayed
support condition.  The GCD sum is the same for every fixed derivative
order.
\end{theorem}

\begin{corollary}\label{cor:central-exponent}
Assume GRH and fix $\ell\geq0$.  Then
\[
 \max_{\substack{X<d\leq2X\\d\in\Fplus}}
 (-1)^\ell L^{(\ell)}(\tfrac12,\chi_d)
 \geq\exp\!\left((1+o(1))
 \sqrt{\frac{\log X\log_3X}{\log_2X}}\right).
\]
\end{corollary}

The next result concerns a shrinking band to the right of the center.
In the resonance calculation, we select square main terms from each
half of the approximate functional equation.  We call the resulting
two terms the \emph{retained contributions}.
A \emph{tilt} is an exponential change of weight. Here the two opposite
tilts replace $p^{-1/2}$ by $p^{-1/2\mp\alpha}$ in the local factors.
We use \emph{free energy} for the logarithm of a normalized finite Euler
product.  A relative error is measured after division by the indicated
positive sum or product.

We use the scales
\[
 \mathcal V_X=\sqrt{\frac{\log X\log_3X}{\log_2X}},\qquad
 \mathcal M_X=\sqrt{\log X\log_2X\log_3X}.
\]
The finite products and their free energies are defined in
Lemma~\ref{lem:mesoscopic-resonator}.

\begin{theorem}\label{thm:mesoscopic}
Assume GRH, fix $A_0>0$ and $\ell\geq0$, and fix
$0<\beta<1/4$, $0<\gamma<1$, and $a>1$ with $a\gamma<1$.
Let $0\leq A\leq A_0$ be allowed to depend on $X$, and put
\[
 N=\lfloor X^{1/4-\beta}\rfloor,\qquad
 \alpha=\frac A{\log_2X},\qquad \sigma_X=\frac12+\alpha.
\]
Let $f,\mathcal M,G_N,G_N^-$, and $\mu_N$ be as in
Lemma~\ref{lem:mesoscopic-resonator}.  Then
\begin{align}\label{eq:mesoscopic-main}
 \max_{\substack{X<d\leq2X\\d\in\Fplus}}
 (-1)^\ell L^{(\ell)}(\sigma_X,\chi_d)
 &\geq(1-o(1))\left\{
 \mu_N^\ell\exp(G_N)\right.\notag\\
 &\hspace{26mm}\left.
 +X^{-\alpha}(\log X)^\ell\exp(G_N^-)\right\}.
\end{align}
The convention is $\mu_N^0=1$.  The relative error is uniform for
$0\leq A\leq A_0$ with the other parameters fixed.
\end{theorem}

In \eqref{eq:mesoscopic-main}, the first contribution carries the
typical logarithmic derivative weight $\mu_N^\ell$.
The dual contribution includes the conductor factor $X^{-\alpha}$.

\begin{corollary}\label{cor:mesoscopic-exponent}
Assume GRH, fix $A_0>0$ and $\ell\geq0$, let $0\leq A\leq A_0$ be
allowed to depend on $X$, and set $\sigma_X=1/2+A/\log_2X$.
Then
\[
 \max_{\substack{X<d\leq2X\\d\in\Fplus}}
 (-1)^\ell L^{(\ell)}(\sigma_X,\chi_d)
 \geq\exp\bigl((\mathrm e^{-A}+o(1))\mathcal V_X\bigr),
\]
uniformly for $0\leq A\leq A_0$.
\end{corollary}

For the distinction between the relative estimate and its
exponential consequence, see Remark~\ref{rem:weight-precision}.
On a smaller boundary scale, the two free energies differ by $o(1)$.
A common finite-product normalization then gives the following result.

\begin{theorem}\label{thm:boundary-layer}\label{prop:boundary-transfer}
Assume GRH and fix $\ell\geq0$, $0<\beta<1/4$, $0<\gamma<1$, and
$a>1$ with $a\gamma<1$.  Let $A\geq0$ satisfy
$A\mathcal V_X\to0$, put
\[
 N=\lfloor X^{1/4-\beta}\rfloor,\qquad
 \alpha=\frac A{\log_2X},\qquad \sigma_X=\frac12+\alpha,
\]
and let $f,\mathcal M,\mu_N$ be as in
Lemma~\ref{lem:mesoscopic-resonator}.
Let $G_N^0$ be as in Corollary~\ref{cor:opposite-tilts}.  Then
\begin{equation}\label{eq:boundary-precise}
 \max_{\substack{X<d\leq2X\\d\in\Fplus}}
 (-1)^\ell L^{(\ell)}(\sigma_X,\chi_d)
 \geq(1-o(1))\exp(G_N^0)
 \left\{\mu_N^\ell+X^{-\alpha}(\log X)^\ell\right\}.
\end{equation}
The convention is $\mu_N^0=1$.  For every prescribed sequence
$a_X\geq0$ with $a_X\mathcal V_X\to0$, the relative error is uniform
for $0\leq A\leq a_X$.
Moreover, on every such range there are admissible choices of
$\beta,\gamma,a$, tending respectively to $0,1,1$, for which
\eqref{eq:boundary-precise} remains uniform and
\[
 \frac{\mu_N}{\mathcal M_X}\longrightarrow1,\qquad
 G_N^0=(1+o(1))\mathcal V_X
\]
uniformly for $0\leq A\leq a_X$.
\end{theorem}

For these diagonal choices and $\ell\geq1$, the logarithm of the
ratio of the dual contribution to the first contribution is
\[
 -\frac{A\log X}{\log_2X}
 +\frac\ell2(\log_2X-\log_3X-\log_4X)+o(1).
\]
The asymptotic crossover location is
\begin{equation}\label{eq:boundary-location}
 A_{\ell,*}(X):=
 \frac{\ell\log_2X}{2\log X}
 \bigl(\log_2X-\log_3X-\log_4X\bigr).
\end{equation}
Within the range $A\mathcal V_X\to0$, consider the quantity
\[
 (A-A_{\ell,*})\frac{\log X}{\log_2X}.
\]
If it tends to $-\infty$, the dual contribution dominates.
If it tends to $+\infty$, the first contribution dominates.
If it remains bounded, the two contributions have the same order.
The common free energy $G_N^0$ makes this comparison possible.
Since $A_{\ell,*}\mathcal V_X\to0$, the entire crossover lies in
the range of Theorem~\ref{thm:boundary-layer}.
These comparisons concern the two terms in the lower bound
\eqref{eq:boundary-precise}.
For $\ell=0$, their formal equality occurs at $A=0$, so there is
no positive crossover location.

\begin{corollary}\label{cor:boundary-exponent}
Assume GRH, fix $\ell\geq0$, and let $A\geq0$ satisfy
$A\mathcal V_X\to0$.  With $\sigma_X=1/2+A/\log_2X$, one has
\begin{equation}\label{eq:boundary-layer}
 \max_{\substack{X<d\leq2X\\d\in\Fplus}}
 (-1)^\ell L^{(\ell)}(\sigma_X,\chi_d)
 \geq\exp\bigl((1+o(1))\mathcal V_X\bigr).
\end{equation}
For every sequence $a_X\geq0$ with $a_X\mathcal V_X\to0$, the error
may be chosen uniformly for $0\leq A\leq a_X$.
\end{corollary}

For the fixed-$\sigma$ result, the truncation and the resonator support
must be chosen together.  For every integer $r\geq1$ with
$\sigma>1-1/r$, put
\[
 \kappa_r(\sigma)=
 \frac{r+1}{4r^2(\sigma-1+1/r)}.
\]
The remaining exponent after the Burgess truncation is
\begin{equation}\label{eq:Delta-r}
 \Delta_r(\sigma):=1-2\kappa_r(\sigma)(1-\sigma)
 =1-\frac{(r+1)(1-\sigma)}
 {2r^2(\sigma-1+1/r)}.
\end{equation}
Let $\Delta(\sigma)$ be the maximum of the positive numbers
$\Delta_r(\sigma)$, and define
\begin{equation}\label{eq:C-Burgess}
 C_{\rm B}(\sigma):=
 \frac1{1-\sigma}
 \left(\frac{\Delta(\sigma)}{2\log2}\right)^{1-\sigma}.
\end{equation}
There are only finitely many eligible $r$ for fixed $\sigma$.
The choice $r=1$ always gives
$\Delta_1(\sigma)=(2\sigma-1)/\sigma>0$, so this maximum exists.

\begin{theorem}\label{thm2}
Assume GRH, and fix an integer $\ell\geq0$ and $1/2<\sigma<1$.
Choose fixed $0<b<1$ and an integer $r\geq1$ with $\sigma>1-1/r$,
and set
\[
 \alpha(b)=2\log\frac{(1+b)^2}{1+b^2}.
\]
Suppose the fixed parameters $\kappa,a$ satisfy
\[
 \kappa>\kappa_r(\sigma),\qquad
 0<a<\frac{1-2\kappa(1-\sigma)}{\alpha(b)},
 \qquad Y=a\log X\log_2X.
\]
For $p\leq Y$, define
\[
 c_p=
 \frac{2b\,h(p)}
 {(1-h(p))(1-b^2)^2+h(p)(1+b^2)},
\]
and put
\[
 Z_Y(\sigma)=\prod_{p\leq Y}(1+c_pp^{-\sigma}),\qquad
 \mu_Y=\sum_{p\leq Y}
 \frac{c_pp^{-\sigma}}{1+c_pp^{-\sigma}}\log p.
\]
Then
\begin{equation}\label{eq:fixed-main-relative}
 \max_{\substack{X<d\leq2X\\d\in\Fplus}}
 (-1)^\ell L^{(\ell)}(\sigma,\chi_d)
 \geq(1-o(1))\mu_Y^\ell Z_Y(\sigma).
\end{equation}
The error tends to zero for each fixed admissible parameter choice.
For $\ell=0$, the convention is $\mu_Y^0=1$.
\end{theorem}

The derivative factor $\mu_Y^\ell$ in \eqref{eq:fixed-main-relative}
is measured against the same finite product $Z_Y(\sigma)$ for every
$\ell$. See Remark~\ref{rem:weight-precision}.
The prime-sum evaluations in the proof give
\[
 \mu_Y=\left(\frac{2b}{1+b^2}\frac1{1-\sigma}+o(1)\right)Y^{1-\sigma},
 \qquad
 \log Z_Y(\sigma)=
 \left(\frac{2b}{1+b^2}\frac1{1-\sigma}+o(1)\right)
 \frac{Y^{1-\sigma}}{\log Y}.
\]
Optimizing the fixed admissible parameters yields the following
coarser exponential consequence.

\begin{corollary}\label{cor:fixed-exponent}
Assume GRH, and let $\ell\geq0$ and $1/2<\sigma<1$ be fixed.  Then
\[
 \max_{\substack{X<d\leq2X\\d\in\Fplus}}
 (-1)^\ell L^{(\ell)}(\sigma,\chi_d)
 \geq
 \exp\!\left((C_{\rm B}(\sigma)+o(1))
 \frac{(\log X)^{1-\sigma}}{(\log_2X)^\sigma}\right).
\]
\end{corollary}

Optimizing the Burgess parameter improves the
P\'olya--Vinogradov truncation whenever an eligible $r>1$ gives a
larger $\Delta_r(\sigma)$ than $r=1$.
For example, at $\sigma=0.9$, the choices $r=1$ and $r=2$ give,
respectively, $\Delta_1=8/9$ and $\Delta_2=29/32$.
The derivative weight in Theorem~\ref{thm2} comes from the typical
value of $\log k$ in the resonating Euler product.
Proposition~\ref{prop:fixed-transfer} proves this statement at the level
of the weighted moment.  Retaining only one prime near the top of the
support would give merely $(\log_2X)^\ell$ in the comparison with the
same product.

\begin{remark}\label{rem:weight-precision}
For positive functions $P(X)$ and $F(X)$ with $\log P(X)=o(F(X))$,
\[
 P(X)\exp\bigl((1+o(1))F(X)\bigr)
 =\exp\bigl((1+o(1))F(X)\bigr).
\]
An estimate of this form cannot detect a separate lower-order
improvement from $P(X)$ when the exponent has an unquantified
$o(F(X))$ error.
For example, fixed powers of $\log X$ are absorbed by the error on the
scale $\mathcal V_X$, and fixed powers of $\log X\log_2X$ are absorbed
by the error on the scale
$(\log X)^{1-\sigma}/(\log_2X)^\sigma$ for fixed $\sigma<1$.

Theorems~\ref{thm1}, \ref{thm:mesoscopic},
\ref{thm:boundary-layer}, and~\ref{thm2} therefore compare derivative
weights with specified GCD sums or finite products before taking their
leading exponential asymptotics.
In the boundary estimate, the common factor $\exp(G_N^0)$ also permits
a comparison of the two contributions.
Corollaries~\ref{cor:central-exponent}, \ref{cor:mesoscopic-exponent},
\ref{cor:boundary-exponent}, and~\ref{cor:fixed-exponent} give only the
leading exponential growth.
The Dickman-transform bounds near the $1$-line keep their stated
relative precision.
\end{remark}

For $A\geq0$ define
\begin{equation}\label{eq:Dell}
 D_\ell(A):=\int_0^\infty \e^{Au}u^\ell\rho(u)\dd u,
\end{equation}
where $\rho$ is the Dickman function.  This integral is finite for every
fixed $A$.

\begin{theorem}\label{thm3}
Assume GRH.  Let $A\geq0$ and $\ell\geq0$ be fixed, and set
$\sigma_A=1-A/\log_2X$.  Then
\[
 \max_{\substack{X<d\leq2X\\d\in\Fplus}}
 (-1)^\ell L^{(\ell)}(\sigma_A,\chi_d)
 \geq\bigl(D_\ell(A)+o(1)\bigr)(\log_2X)^{\ell+1}.
\]
\end{theorem}

\begin{corollary}\label{thm4}
Assume GRH and let $\ell\geq0$ be fixed.  Then
\[
 \max_{\substack{X<d\leq2X\\d\in\Fplus}}
 (-1)^\ell L^{(\ell)}(1,\chi_d)
 \geq\bigl(D_\ell(0)+o(1)\bigr)(\log_2X)^{\ell+1}.
\]
\end{corollary}

\begin{remark}\label{rem:log-derivative}
For $\ell=1$ the identity
$-L'(s,\chi)=L(s,\chi)(-L'/L)(s,\chi)$ suggests a possible refinement of
the lower-order terms at the $1$-line.
Recent sharp resonance bounds for the logarithmic derivative appear in
\cite{LogDerivatives}.  Separate maximizers for the two factors do not
give a bound for their product.
Refining the lower-order terms in Corollary~\ref{thm4} would require a
joint resonance moment.  This remains a direction for further work.
\end{remark}

The transition result is most naturally stated using
\[
 A=(1-\sigma)\log_2X,
 \qquad B=(\sigma-\tfrac12)\log_2X.
\]
In the first part below, $D_\ell(A)$ is allowed to depend on $X$ through
$A$.

\begin{theorem}\label{thm5}
Assume GRH, let $\ell\geq0$ be fixed, and suppose
\[
 \frac12<\sigma<1-\frac1{\log_2X},\qquad B\longrightarrow\infty.
\]
If $1\leq A\leq\frac12\log_3X$, then
\begin{equation}\label{eq:transition-strong}
 \max_{\substack{X<d\leq2X\\d\in\Fplus}}
 (-1)^\ell L^{(\ell)}(\sigma,\chi_d)
 \geq (1-o(1))D_\ell(A)(\log_2X)^{\ell+1}.
\end{equation}
For the remaining range $A>\frac12\log_3X$, one has
\begin{equation}\label{eq:transition-standard}
 \max_{\substack{X<d\leq2X\\d\in\Fplus}}
 (-1)^\ell L^{(\ell)}(\sigma,\chi_d)
 \gg_\ell
 \exp\!\left(\frac{\e^A}{A}-A\right)
 (\log_2X)^{\ell+1}.
\end{equation}
\end{theorem}

The first conclusion is stronger than \eqref{eq:transition-standard}.
For instance, the Dickman transform satisfies
\begin{equation}\label{eq:D0-transform}
 D_0(A)=\exp\!\left(\gamma+
 \int_0^A\frac{\e^t-1}{t}\dd t\right)
 =\frac{\exp(\operatorname{Ei}(A))}{A}\qquad(A>0).
\end{equation}
The first identity is the Dickman Laplace-transform formula
\cite[Lemma~2.6]{HT}. The second is its standard exponential-integral
rewriting.

\begin{remark}\label{rem:central-band}
The condition $B\to\infty$ in Theorem~\ref{thm5} is essential in its
proof.  The bounded part of the excluded band,
$\sigma=1/2+A/\log_2X$ with $A=O(1)$, is covered by the relative
transfer in Theorem~\ref{thm:mesoscopic} and its exponential
Corollary~\ref{cor:mesoscopic-exponent}, including $A=0$.
On the smaller boundary scale, Theorem~\ref{thm:boundary-layer}
compares the two derivative contributions through a common free
energy.  The opposite tilts enter the comparison.
For the distinction from the exponential consequence, see
Remark~\ref{rem:weight-precision}.
\end{remark}

\medskip
\noindent\textbf{Comparison with earlier and recent results.}
The following comparisons concern the leading scale.
The underlying families and uniformity ranges need not be identical.

At the central point, Soundararajan's resonance method gave the earlier
bound
\[
 \exp\left(\left(\frac1{\sqrt5}+o(1)\right)
 \sqrt{\frac{\log X}{\log_2X}}\right)
\]
for quadratic characters \cite{SoundRes}.  Darbar and Maiti introduced
the GRH-based quadratic long-resonator transfer and obtained the stronger
$\sqrt{\log_3X}$ scale with constant $1/2$ in the exponent
\cite{DarbarMaiti}.  The recent preprint \cite{ExtremeQuadratic} raises
that constant to $1$.  Corollary~\ref{cor:central-exponent} gives this
leading exponential scale for every fixed derivative order in the
positive even subfamily $\Fplus$.
Theorem~\ref{thm1} compares each fixed derivative with the same GCD sum,
including the factor $(\log X)^\ell$. See
Remark~\ref{rem:weight-precision}.
A subsequent preprint treats quadratic values with prime-related moduli
\cite{PrimeRelated}.  Its arithmetic setting uses a sparser set of
moduli and does not address derivatives.

For fixed $1/2<\sigma<1$ and $\ell=0$, Darbar--Maiti's logarithmic Euler
product argument gives the constant
\begin{equation}\label{eq:CDM}
 C_{\rm DM}(\sigma)=\frac1{1-\sigma}
 \left(\frac1{2\log 2}\right)^{1-\sigma}
\end{equation}
for the full quadratic family \cite{DarbarMaiti}.
Lamzouri had earlier obtained related distributional results and random
Euler-product predictions \cite{Lamzouri}.
Since $\Delta(\sigma)<1$, one has
$C_{\rm B}(\sigma)<C_{\rm DM}(\sigma)$.
For $\ell=0$, the constant in Corollary~\ref{cor:fixed-exponent}
is therefore smaller than this known constant.
Theorem~\ref{thm2} gives a relative lower bound for every fixed
derivative against the finite product $Z_Y(\sigma)$, with the full
typical factor $\mu_Y^\ell$.
For fixed admissible parameters this mean has order
$(\log X\log_2X)^{1-\sigma}$.
Within this differentiated
Dirichlet-series method, optimizing the Burgess parameter strictly
improves the P\'olya--Vinogradov choice whenever an eligible $r>1$
gives $\Delta_r(\sigma)>\Delta_1(\sigma)$.  For comparison,
\cite[Theorem~1.1]{CuiPengSongZhao} gives exponential lower bounds on this
scale for fixed derivatives when one maximizes over all nonprincipal
characters modulo a prime.  Our result isolates quadratic characters and records the typical
derivative weight relative to a fixed finite Euler product.

Near the $1$-line, Yang proved Dickman-constant lower bounds for
derivatives when one maximizes over all nonprincipal characters modulo
$q$ \cite{YangDerivatives}. The systematic transition beyond the
$1/\log_2q$ scale was developed further in \cite{DSWZ}.  Theorems
\ref{thm3}--\ref{thm5} recover the same leading Dickman constants inside
the sparser family of real quadratic characters, under GRH.  The first part of
Theorem~\ref{thm5} is stronger than the previously available general
transition lower bound because it keeps the full transform $D_\ell(A)$
uniformly for the growing range $A\leq\frac12\log_3X$.  The second part
matches the established transition scale but is uniform down to
$B\to\infty$.

In the shrinking central band, Theorem~1.1 of \cite{NearCenter} states a
value bound ($\ell=0$) with exponent constant
$1/(2(\e-1)\e^A)$.  Our two-sided construction gives
$\mathrm e^{-A}$ in Corollary~\ref{cor:mesoscopic-exponent} for every
fixed derivative.  Theorem~\ref{thm:mesoscopic} includes the dual half.
Theorem~\ref{thm:boundary-layer} locates the crossover between the
two derivative contributions at \eqref{eq:boundary-location}.
Theorem~\ref{thm5} treats the complementary regime $B\to\infty$.

\begin{remark}\label{rem:sharpness}
The central exponential constant in
Corollary~\ref{cor:central-exponent} agrees with the best currently
available quadratic result, and the mesoscopic constant
$\mathrm e^{-A}$ in Corollary~\ref{cor:mesoscopic-exponent} joins it
continuously as $A\downarrow0$.
Theorem~\ref{thm:boundary-layer} compares the polynomial derivative
weights in this lower bound.
Determining the asymptotics of the maximum would also require matching
upper bounds, which remain open in the central and mesoscopic bands.
Here the friable character-sum model predicts maximal character sums
by smooth-number counts, whose limiting profile is the Dickman function.
At the $1$-line and on its fixed-$A$ scale, the constants
$D_\ell(A)$ agree with those predicted by this model.
The conjectural sharp upper bounds would follow from the uniform
character-sum approximation proposed by Granville and Soundararajan
\cite{GSchar}.
For fixed $\sigma$, the comparison \eqref{eq:CDM} concerns values
without derivatives.  Optimizing bounds for derivatives is a separate
problem.  Another open question is a sharp joint resonance estimate for
$L$ and its logarithmic derivatives; see
Remark~\ref{rem:log-derivative} and \cite{LogDerivatives}.
\end{remark}

\section{Preliminary lemmas}\label{sec2}

For a positive integer $n=n_0n_1^2$, with $n_0$ square-free, set
\[
 g_1(n_0)=\exp\bigl((\log(2n_0))^{1-\eps}\bigr),\qquad
 g_2(n_1)=\sum_{a\mid n_1}\frac{\mu(a)^2}{a^{1/2+\eps}}.
\]

\begin{lemma}\label{lem:orth}
Assume GRH.  Uniformly for $x\geq2$ and $n\geq1$,
\begin{equation}\label{eq:orth}
 \sum_{\substack{d\leq x\\d\in\Fplus}}\chi_d(n)
 =\frac{2x}{\pi^2}h(n)\ind_{n=\square}
 +O_\eps\!\left(x^{1/2+\eps}g_1(n_0)g_2(n_1)\right).
\end{equation}
The same estimate holds over $X<d\leq2X$ by subtraction.
\end{lemma}

Lemma~\ref{lem:orth} is the congruence-class form of \cite[Lemma~1]{DarbarMaiti}.
Although the cited lemma is stated for the full family of fundamental
discriminants, its proof already begins with the square-free
$d\equiv1\pmod4$ projection.  We use its character-sum estimate and
record the normalization for the present positive subfamily.
For $j=0,1$, the projected real characters are
\[
 \lambda_j(r)=\chi_{0,4}(r)\chi_4(r)^j
                    \left(\frac r{n_0}\right),
\]
where $\chi_{0,4}$ and $\chi_4$ are the principal and
nonprincipal characters modulo $4$, respectively.  The symbol is the
Kronecker symbol.
Their primitive conductors divide
$8n_0$, and both are nonprincipal when $n_0>1$.
When $2\mid n_0$, the $2$-part is $\chi_8$ or $\chi_{-8}$.
The same GRH estimate in the cited proof applies.  Its
square-part divisor majorant is $g_2(n_1)$.
The fixed modulus factors are absorbed by $g_1(n_0)$ using the same
slack in the exponent.
When $n_0=1$, only $\lambda_0=\chi_{0,4}$ is principal.
The density of square-free positive integers in the class $1\pmod4$
is
\[
 \frac14\prod_{p>2}(1-p^{-2})=\frac2{\pi^2}.
\]
For square $n$, requiring $(d,n)=1$ contributes $p/(p+1)$ at each
odd prime dividing $n$.  There is no factor at $2$, since every
$d\in\Fplus$ is odd.  This gives the main term in
\eqref{eq:orth}.  Subtraction gives the interval form.

We next improve the truncation used in the fixed-$\sigma$ range.

\begin{lemma}\label{lem:burgess-trunc}
Fix $1/2<\sigma<1$, an integer $r\geq1$ with
$\sigma>1-1/r$, and an integer $\ell\geq0$.  Put
\[
 \kappa_r(\sigma):=
 \frac{r+1}{4r^2(\sigma-1+1/r)}.
\]
For every fixed $\kappa>\kappa_r(\sigma)$,
\begin{equation}\label{eq:burgess-trunc}
 (-1)^\ell L^{(\ell)}(\sigma,\chi_d)
 =\sum_{n\leq d^\kappa}
 \frac{\chi_d(n)(\log n)^\ell}{n^\sigma}+o(1),
\end{equation}
uniformly for $X<d\leq2X$, $d\in\Fplus$.  For $r=1$,
$\kappa_1(\sigma)=1/(2\sigma)$.
\end{lemma}

For $r\geq2$, this follows from Burgess's character-sum estimate
\cite{Burgess} by partial summation. The estimate applies because $d$ is
square-free.  For $r=1$, it follows in the same way from the
P\'olya--Vinogradov inequality \cite[Chapter~12]{IK}.
With $K=d^\kappa$, the corresponding tail bounds are
\[
 \sum_{n>K}\frac{\chi_d(n)(\log n)^\ell}{n^\sigma}
 \ll
 \begin{cases}
 d^{(r+1)/(4r^2)+\eps}
 K^{1-1/r-\sigma}(\log(2K))^\ell,&r\geq2,\\
 d^{1/2}K^{-\sigma}(\log d)(\log(2K))^\ell,&r=1.
 \end{cases}
\]
The Burgess bound for initial character sums extends beyond one period
because the sum over a complete period is zero.  Taking a sufficiently
small fixed $\eps>0$ gives the asserted uniform $o(1)$.

\begin{remark}\label{rem:burgess-parameter}
The integer $r$ does not index a hierarchy in which larger is always
better.  The condition $\sigma>1-1/r$ ensures that the Burgess tail
decays after partial summation.  The resonator moment also requires
$1-2\kappa_r(\sigma)(1-\sigma)>0$.  These two constraints leave only
finitely many useful choices for a fixed $\sigma$.  The quantity
$\Delta_r(\sigma)$ in \eqref{eq:Delta-r} is precisely the remaining
exponent budget after truncation.  The choice $r=1$ gives the
P\'olya--Vinogradov baseline.  The maximum defining
$\Delta(\sigma)$ automatically selects a genuine Burgess improvement
when one is available.
\end{remark}

We quote \cite[Lemma~6]{BondarenkoSeipNote} in the following form.

\begin{lemma}\label{lem:prime-sum}
There is an absolute constant $C$ such that, uniformly for
$1/2<\sigma<1$ and $(1-\sigma)\log y\geq1/2$,
\begin{equation}\label{eq:prime-sum}
 \sum_{p\leq y}\frac1{p^\sigma}
 \geq \sigma\log_2y+
 \frac{y^{1-\sigma}}{(1-\sigma)\log y}-C.
\end{equation}
\end{lemma}

We express the derivative weights through derivatives of a finite
Euler product.  The following elementary lemma gives the form used below.

\begin{lemma}\label{lem:free-energy}
Fix $1/2<\sigma<1$, $0<c\leq1$, and an integer $\ell\geq0$.  Suppose that
$c_p\geq0$ and $c_p=c+O(p^{-1})$ for $p\leq Y$, and set
\[
 Z_Y(\sigma)=\prod_{p\leq Y}(1+c_pp^{-\sigma}),\qquad
 \mu_Y=\sum_{p\leq Y}
 \frac{c_pp^{-\sigma}}{1+c_pp^{-\sigma}}\log p.
\]
Then
\begin{align}
 \sum_{\substack{k\ {\rm square\mbox{-}free}\\P^+(k)\leq Y}}
 \frac{c_k(\log k)^\ell}{k^\sigma}
 &=(1+o(1))\mu_Y^\ell Z_Y(\sigma),\label{eq:free-energy}\\
 \mu_Y&=\left(\frac{c}{1-\sigma}+o(1)\right)Y^{1-\sigma},
 \label{eq:mu-asymp}
\end{align}
where $c_k=\prod_{p\mid k}c_p$.  For $\ell=0$, the first factor on the
right of \eqref{eq:free-energy} is $1$.
\end{lemma}

\begin{proof}
Give each prime $p\leq Y$ an independent Bernoulli variable $I_p$ with
\[
 \mathbb P(I_p=1)=\frac{c_pp^{-\sigma}}{1+c_pp^{-\sigma}}.
\]
After division by $Z_Y(\sigma)$, the left side of
\eqref{eq:free-energy} is the $\ell$-th moment of
$S_Y=\sum_{p\leq Y}I_p\log p$.  Its mean is $\mu_Y$, and
\[
 \operatorname {Var}(S_Y)
 \leq\sum_{p\leq Y}\frac{c_p(\log p)^2}{p^\sigma}
 \ll_{\sigma,c}Y^{1-\sigma}\log Y.
\]
For the mean, expanding the local denominator gives
\[
 \mu_Y=c\sum_{p\leq Y}\frac{\log p}{p^\sigma}
 +O\!\left(\sum_p\frac{\log p}{p^{\sigma+1}}
            +\sum_p\frac{\log p}{p^{2\sigma}}\right).
\]
Both error sums converge because $\sigma>1/2$.
Partial summation of the prime number theorem gives
$\sum_{p\leq Y}(\log p)p^{-\sigma}
=(1+o(1))Y^{1-\sigma}/(1-\sigma)$.
This proves \eqref{eq:mu-asymp} and gives
$\operatorname {Var}(S_Y)=o(\mu_Y^2)$.  For every fixed $j\geq2$, the
$j$-th cumulant of $S_Y-\mu_Y$ is
\[
 \ll_j\sum_{p\leq Y}\frac{(\log p)^j}{p^\sigma}
 \ll_{j,\sigma}Y^{1-\sigma}(\log Y)^{j-1}.
\]
For each fixed $j\geq2$, division by $\mu_Y^j$ makes this bound
$O_{j,\sigma}( (\log Y/Y^{1-\sigma})^{j-1})=o(1)$.
The centered variable has first cumulant zero.  In the formula for its
even moments as sums over partitions, every nonzero term is
a product of cumulants of orders at least two.  After normalization by
the corresponding power of $\mu_Y$, every such product tends to zero.
There are only finitely many partitions for each fixed order, so
$\mathbb E|S_Y-\mu_Y|^{2j}=o(\mu_Y^{2j})$ for every fixed $j\geq1$.
H\"older's inequality gives
$\mathbb E|S_Y-\mu_Y|^j=o(\mu_Y^j)$ for every fixed $j\geq1$.
Expanding
\[
 \mathbb E S_Y^\ell
 =\mu_Y^\ell+
 \sum_{j=1}^{\ell}\binom\ell j\mu_Y^{\ell-j}
             \mathbb E(S_Y-\mu_Y)^j
\]
shows that every term in the sum is $o(\mu_Y^\ell)$.
Multiplication by $Z_Y(\sigma)$ proves \eqref{eq:free-energy}.
\end{proof}

For later use, write $\Ssmooth(y)=\{n\geq1:P^+(n)\leq y\}$.

\begin{lemma}\label{lem:friable}
Let $\ell\geq0$ be fixed, $y\to\infty$, and
$\sigma=1-A/\log y$.  Uniformly for $A$ in a fixed compact subset of
$[0,\infty)$,
\begin{equation}\label{eq:friable-fixed}
 \sum_{n\in\Ssmooth(y)}\frac{(\log n)^\ell}{n^\sigma}
 =\bigl(D_\ell(A)+o(1)\bigr)(\log y)^{\ell+1}.
\end{equation}
The same formula, with a relative $o(1)$ error, remains valid uniformly for
$1\leq A\leq\frac12\log_3X$ when
\[
 y=\frac{\log X}{100\log_2X},\qquad
 \sigma=1-\frac A{\log_2X}.
\]
For $y=\log X/(100\log_2X)$ and
$\sigma=1-A/\log_2X$, the fixed-$A$ and growing-$A$ assertions are
unchanged if the sum is restricted to $n\leq X^c$ for any fixed $c>0$.
\end{lemma}

\begin{proof}
The fixed-parameter weighted Dickman calculation is given in
\cite[proof of Proposition~1, Section~5.1 of the preprint]{DSWZ}.
Its specialization at $A=0$ appears in
\cite[proof of Proposition~1, Section~3 of the preprint]
{YangDerivatives}.  We give the estimates needed to extend this calculation
to the full series, obtain uniformity in the parameter, and justify truncation.
Write $L=\log y$ and keep $(1-\sigma)L$ explicit.  By
\cite[Theorem~1.1 and Corollary~2.3]{HT},
\begin{equation}\label{eq:HT-uniform}
 \Psi(y^u,y)=y^u\rho(u)
 \left(1+O\left(\frac{\log(u+1)}L\right)\right)
 \qquad(1\leq u\leq L^{2/3}),
\end{equation}
and
\begin{equation}\label{eq:Dickman-tail}
 \log\rho(u)=-u\left(\log u+\log_2(u+2)-1
 +O\left(\frac{\log_2(u+2)}{\log(u+2)}\right)\right)
 \qquad(u\geq1).
\end{equation}
Outside this range, we also need an upper bound for the friable count:
\begin{equation}\label{eq:Psi-tail}
 \frac{\Psi(y^u,y)}{y^u}\ll
 \begin{cases}
  \exp(-9u\log u/10),
          &L^{2/3}\leq u\leq\sqrt y,\\
  \exp(-uL/3),&u>\sqrt y.
 \end{cases}
\end{equation}
The first bound follows directly from \cite[Corollary~1.3]{HT},
with a fixed $\varepsilon=1/4$, since $u\geq L^{2/3}\to\infty$
and $u\leq\sqrt y\leq y^{1-\varepsilon}$.  In the second range,
Rankin's inequality with exponent $1/2$ gives
$-uL/2+O(\sqrt y/\log y)\leq-uL/3$, proving \eqref{eq:Psi-tail}.

Partial summation of the weighted series gives
\begin{equation}\label{eq:friable-Stieltjes}
 \sum_{n\in\Ssmooth(y)}\frac{(\log n)^\ell}{n^\sigma}
 =L^{\ell+1}\int_0^\infty
 \frac{\Psi(y^u,y)}{y^u}e^{(1-\sigma)Lu}
 \left(\sigma u^\ell-\frac{\ell u^{\ell-1}}L\right)\dd u,
\end{equation}
where the second term is omitted for $\ell=0$.
The boundary term at infinity vanishes by Rankin's inequality with
any positive exponent below $\sigma$.  Equations
\eqref{eq:HT-uniform}--\eqref{eq:Psi-tail}, together with
$\Psi(y^u,y)/y^u=1+O(y^{-u})$ on $0\leq u\leq1$, give a uniform version
of the cited fixed-parameter calculation on compact sets of $(1-\sigma)L\geq0$.
This proves \eqref{eq:friable-fixed}.  Here $\sigma=1+O(L^{-1})$,
$D_{\ell-1}((1-\sigma)L)\ll_\ell D_\ell((1-\sigma)L)$ for $\ell\geq1$, and the contribution
of $O(y^{-u})$ is $O(L^{-1})$ after removal of $L^{\ell+1}$.

We now prove the growing-parameter assertion.  For the specified $y$,
\[
 L=\log_2X-\log_3X+O(1),\qquad
 (1-\sigma)L=A+O\left(\frac{A\log_3X}{\log_2X}\right).
\]
The exact Laplace transform \cite[Lemma~2.6]{HT} is
\[
 D_0(v)=\exp\left(\gamma+\int_0^v\frac{e^t-1}{t}\dd t\right).
\]
For each fixed integer $j\geq1$, its logarithmic derivatives satisfy
\[
 0\leq\frac{\mathrm d^j}{\mathrm dv^j}\log D_0(v)
 =\int_0^1t^{j-1}e^{vt}\dd t
 \leq\frac{e^v-1}{v}\qquad(v\geq0).
\]
The ratio $(e^v-1)/v$ is interpreted as $1$ at $v=0$ and is at least $1$.
The nonnegative Bell polynomial for $D_j(v)/D_0(v)$ has leading term
$((e^v-1)/v)^j$.  Every other term is at most a fixed multiple of it.
We obtain, also for $j=0$,
\begin{equation}\label{eq:D-ratio}
 \frac{D_{j+1}(v)}{D_j(v)}\ll_j\frac{e^v-1}{v}
 \ll_j1+\frac{e^v}{1+v}\qquad(v\geq0).
\end{equation}
The mean-value theorem for $\log D_\ell$ gives, uniformly for
$1\leq A\leq\tfrac12\log_3X$,
\[
 D_\ell((1-\sigma)L)=D_\ell(A)(1+o(1)),\qquad
 |(1-\sigma)L-A|\left(1+\frac{e^A}{1+A}\right)
 \ll\frac{\log_3X}{\sqrt{\log_2X}}=o(1).
\]
On $u\leq L^{2/3}$, the relative error in \eqref{eq:HT-uniform} is
$O(\log L/L)=o(1)$.  Beyond this point,
\eqref{eq:Dickman-tail} and \eqref{eq:Psi-tail} bound both
$e^{(1-\sigma)Lu}\rho(u)$ and $e^{(1-\sigma)Lu}\Psi(y^u,y)/y^u$ by
$e^{-c_1u\log_3X}$ up to $\sqrt y$ and by $e^{-c_1uL}$ afterwards.
For $u\geq L^{2/3}$, we have
$\tfrac9{10}\log u-(1-\sigma)L\geq(\tfrac1{10}+o(1))\log L$.
Their integrals with the fixed polynomial weights are $o(D_\ell(A))$.
Also $\sigma=1+o(1)$ and
$D_{\ell-1}((1-\sigma)L)/L=o(D_\ell(A))$ when $\ell\geq1$.
Substitution in \eqref{eq:friable-Stieltjes} proves the growing assertion
with relative error $o(1)$.  The same estimates cover compact $A$.

Finally $c\log X/L>\sqrt y$ for every fixed $c>0$ and large $X$.
Partial summation of the series over $n>X^c$, dropping its nonpositive
boundary term, gives
\[
 \sum_{\substack{n>X^c\\n\in\Ssmooth(y)}}
        \frac{(\log n)^\ell}{n^\sigma}
 \leq\sigma L^{\ell+1}\int_{c\log X/L}^{\infty}
        \frac{\Psi(y^u,y)}{y^u}e^{(1-\sigma)Lu}u^\ell\dd u.
\]
By \eqref{eq:Psi-tail} this is
$O_{c,\ell}(X^{-c/4}(\log X)^\ell)$, uniformly for
$0\leq A\leq\tfrac12\log_3X$.
Since $D_\ell(A)\geq D_\ell(0)>0$, this bound is
$o(D_\ell(A)L^{\ell+1})$.  This proves the truncation assertion.
\end{proof}

\section{Preparations for the proof of Theorem \ref{thm1}}\label{sec3}

For $d\in\Fplus$ and $|\alpha|<1/2$, set
\begin{align*}
 V_\alpha(x)&=\frac1{2\pi i}\int_{(2)}\pi^{-u/2}
 \frac{\Gamma(1/4+\alpha/2+u/2)}{\Gamma(1/4+\alpha/2)}
 x^{-u}\frac{\dd u}{u},\\
 X_\alpha(d)&=(d/\pi)^{-\alpha}
 \frac{\Gamma(1/4-\alpha/2)}{\Gamma(1/4+\alpha/2)}.
\end{align*}

\begin{lemma}
\label{lem:central-afe}
For every fixed $\ell\geq0$ and every $d\in\Fplus$ with $d>1$,
\begin{equation}\label{eq:central-afe}
 (-1)^\ell L^{(\ell)}\!\left(\frac12,\chi_d\right)
 =\sum_{n\geq1}\frac{\chi_d(n)}{\sqrt n}W_{\ell,d}(n),
\end{equation}
where
\begin{equation}\label{eq:Wdef}
 W_{\ell,d}(n)=(-1)^\ell
 \left.\partial_\alpha^\ell\left\{
 n^{-\alpha}V_\alpha(n/\sqrt d)
 +X_\alpha(d)n^\alpha V_{-\alpha}(n/\sqrt d)
 \right\}\right|_{\alpha=0}.
\end{equation}
For every fixed $0<\delta<1/2$,
\begin{equation}\label{eq:Wsmall}
 W_{\ell,d}(n)\gg_{\ell,\delta}(\log X)^\ell
 \qquad(1\leq n\leq X^\delta,\ X<d\leq2X).
\end{equation}
For sufficiently large $X$, one also has $W_{\ell,d}(n)\geq0$ for all
$n\geq1$, and, for every $C>0$,
\begin{equation}\label{eq:Wglobal}
 |W_{\ell,d}(n)|+d\left|\partial_dW_{\ell,d}(n)\right|
 \ll_{\ell,C}(\log(2dn))^\ell
 \left(1+\frac n{\sqrt d}\right)^{-C}.
\end{equation}
\end{lemma}

\begin{proof}
At $\alpha=0$, the definitions in Lemmas~\ref{lem:shifted-kernel}
and \ref{lem:shifted-positive}, proved below, give
\[
 W_{\ell,d}(n)=W^+_{\ell,0,d}(n)+W^-_{\ell,0,d}(n).
\]
Equation~\eqref{eq:shifted-afe} gives
\eqref{eq:central-afe}.

The first-half bound \eqref{eq:shifted-first-global} and the
sharper dual bound for $0\leq\alpha\leq1/4$ in
Lemma~\ref{lem:shifted-kernel}, evaluated at $\alpha=0$, give
\eqref{eq:Wglobal}, including the $d$-derivative estimate.
Here $\log(2Xn)\asymp\log(2dn)$ because $X<d\leq2X$.

Fix $0<\delta<1/2$ and $n\leq X^\delta$.
We have $\pi n^2/d=o(1)$ uniformly.  The interval $1\leq t\leq2$
therefore lies in the dual integral
\eqref{eq:Laplace-two-halves}.  On that interval its polynomial
argument satisfies
\[
 \log\frac d{\pi n}+\frac12\log t
 \geq(1-\delta)\log X+O(1).
\]
The polynomial lower bound in the proof of
Lemma~\ref{lem:shifted-positive} gives a contribution
$\gg_{\ell,\delta}(\log X)^\ell$ from this interval at $\alpha=0$.
The rest of the dual integrand is nonnegative by the same proof.
For the first half, \eqref{eq:shifted-first-small} gives
\[
 W^+_{\ell,0,d}(n)
 =(\log n)^\ell+
 O_{\ell,\delta}(X^{-c_\delta}(\log X)^\ell).
\]
Its main term is nonnegative and its error is negligible relative
to the dual lower bound.  This proves \eqref{eq:Wsmall}.
For $\ell=0$, the polynomial is exactly $1$, so the same argument
applies with $(\log X)^0=1$.

For $n\geq R_\ell$, the two weights at $\alpha=0$ are both
nonnegative by \eqref{eq:dual-positive} and
\eqref{eq:first-positive}.  If $\ell\geq1$ and $1\leq n<R_\ell$,
\eqref{eq:dual-boundary} gives
$W^-_{\ell,0,d}(n)=(1+o(1))(\log X)^\ell$ uniformly over this
finite set.  Equation~\eqref{eq:first-finite-negative} bounds the
first weight below by $-O_\ell(1)$.  Their sum is
positive for large $X$.

Finally, for $\ell=0$, the original definition gives
$W_{0,d}(n)=2V_0(n/\sqrt d)$.  The incomplete-gamma integral
in Lemma~\ref{lem:shifted-kernel} is nonnegative for every $n$.
This completes the proof of the global sign assertion.
\end{proof}

We also need the sharp square-free GCD sum.  There are sets $\M$ of $N$
square-free integers for which
\begin{equation}\label{eq:GCD}
 \sum_{m,n\in\M}\sqrt{\frac{(m,n)}{[m,n]}}
 \geq N\exp\!\left((2+o(1))
 \sqrt{\frac{\log N\log_3N}{\log_2N}}\right),
\end{equation}
and which satisfy
\begin{equation}\label{eq:pmax}
 \max_{m\in\M}P^+(m)\leq(\log N)^{1+o(1)}.
\end{equation}
This is \cite[Eq.~(1.5)]{BT}.  The square-free construction is due to
\cite{BondarenkoSeip}.  Its primes lie in blocks ending at
$\exp((\log_2N)^{\gamma})\log N\log_2N$ with fixed $\gamma<1$, which
gives \eqref{eq:pmax}.  The same support property and the tail estimate
needed below are recorded explicitly in \cite[Section~3]{ExtremeQuadratic}.

\section{Proof of Theorem \ref{thm1}}\label{sec4}

\begin{proof}[Proof of Theorem~\ref{thm1}]
Fix $0<\beta<1/4$ and set
$N=\lfloor X^{1/4-\beta}\rfloor$.
Let $\M$ be any family satisfying the hypotheses of the theorem.
Put
\[
 y_\M=\max_{m\in\M}P^+(m),\qquad
 Q=\prod_{p\leq y_\M}p,
\]
and define
\[
 R_d=\sum_{m\in\M}\chi_d(Qm).
\]
Each $m\in\M$ is square-free and divides $Q$.
Every square diagonal appearing below has precisely the odd
prime divisors of $Q$, and its local density is $h(Q)$.
The prime-support bound implies
\[
 \log Q\ll y_\M\leq(\log N)^{1+o(1)},\qquad
 h(Q)^{-1}=\prod_{\substack{p\leq y_\M\\p>2}}(1+1/p)
 =X^{o(1)}.
\]
These estimates allow the conductor and local-density factors in
Lemma~\ref{lem:orth} to be absorbed into powers $X^{o(1)}$.

Let
\begin{align*}
 S_1&=\sum_{\substack{X<d\leq2X\\d\in\Fplus}}R_d^2,\\
 S_{2,\ell}&=\sum_{\substack{X<d\leq2X\\d\in\Fplus}}
 (-1)^\ell L^{(\ell)}(\tfrac12,\chi_d)R_d^2.
\end{align*}
Expanding $R_d^2$ and applying Lemma~\ref{lem:orth} gives
\begin{equation}\label{eq:central-den}
 S_1=\left(\frac2{\pi^2}+o(1)\right)Xh(Q)N.
\end{equation}
The integer $Q^2mn$ is a square
precisely when $m=n$.  This uses the square-free assumption on $m$ and $n$.
For these $N$ pairs, $h(Q^2m^2)=h(Q)$.
For every pair, the square-free part of $Q^2mn$ divides $Q$.
All prime factors of its square part are at most $y_\M$.
The definitions of $g_1$ and $g_2$ give
\[
 g_1((Q^2mn)_0)g_2((Q^2mn)_1)=X^{o(1)}.
\]
For example, $\log g_1\ll(\log Q)^{1-\eps}=o(\log X)$, and
$\log g_2\leq\sum_{p\leq y_\M}p^{-1/2-\eps}=o(\log X)$.
The total error over the $N^2$ pairs is
$X^{1/2+\eps+o(1)}N^2=o(Xh(Q)N)$ for sufficiently small fixed
$\eps$.

Insert Lemma~\ref{lem:central-afe} into $S_{2,\ell}$ and expand
$R_d^2$.  For each $m,n$ and each positive integer $k$, apply
Lemma~\ref{lem:orth} to $\chi_d(kQ^2mn)$, with the weight
$W_{\ell,d}(k)$.
Partial summation is applicable because \eqref{eq:Wglobal}
controls both this weight and its $d$-derivative.
For any prescribed $\eps>0$, using a smaller positive parameter
in Lemma~\ref{lem:orth} gives, uniformly in $m,n$ and $k$,
\[
 g_1((kQ^2mn)_0)g_2((kQ^2mn)_1)
 \ll_\eps k^\eps X^{o(1)}.
\]
The primes from $Qmn$ contribute $X^{o(1)}$ as above.
The remaining factors of $g_2$ are bounded by an arbitrarily small
power of $k$.  We also have
$(\log(2k))^{1-\eps'}=o(\log k)$ for each fixed $\eps'>0$.
The total orthogonality remainder is bounded
by $N^2X^{1/2+\eps'+o(1)}$ times a sum of the form
\[
 \sum_{k\geq1}k^{-1/2+\eps'}
    \log^\ell(2Xk)(1+k/\sqrt X)^{-C}
 \ll_{\ell,\eps'}X^{1/4+\eps'/2}(\log X)^\ell,
\]
where $\eps'>0$ may be taken arbitrarily small.
It includes the exponent loss in Lemma~\ref{lem:orth}.
We choose $C$ sufficiently large.
To prove this inequality, split at $\sqrt X$.
Above that point, dyadic intervals give a convergent geometric series.
Choosing all intermediate positive exponents sufficiently small
in terms of the displayed $\eps$, we obtain
\begin{equation}\label{eq:central-error}
 \ll_{\ell,\eps}(\log X)^\ell X^{3/4+\eps}N^{2+\eps}.
\end{equation}
This bound includes the orthogonality errors from square as well as
nonsquare arguments.

The square main terms are nonnegative by
Lemma~\ref{lem:central-afe}.  For each $m,n\in\M$, select the integer
\[
 k=\frac{[m,n]}{(m,n)}.
\]
Then
\[
 kQ^2mn=(Q[m,n])^2,\qquad
 h(kQ^2mn)=h(Q),\qquad
 k^{-1/2}=\sqrt{\frac{(m,n)}{[m,n]}}.
\]
We next show that the pairs with $k>X^{1/4}$ make a negligible
contribution.  For fixed $m$, summing over all square-free integers supported
on the primes at most $y_\M$ only increases the sum over $n\in\M$.
At each prime, the two choices for its occurrence in $n$ contribute
$1$ and $p^{-1/3}$, respectively.  These cases correspond to agreement
or disagreement with its occurrence in $m$.
The divisor-product estimate in
\cite[p.~25]{BT} gives
\[
 \sum_{n\in\M}\left(\frac{(m,n)}{[m,n]}\right)^{1/3}
 \leq\prod_{p\leq y_\M}(1+p^{-1/3})
 \leq\exp\!\bigl(y_\M^{2/3+o(1)}\bigr)=N^{o(1)}.
\]
For $k>X^{1/4}$ we have
$k^{-1/2}=k^{-1/3}k^{-1/6}\leq X^{-1/24}k^{-1/3}$.
Summing this inequality over $m$ proves
\[
 \sum_{\substack{m,n\in\M\\k(m,n)>X^{1/4}}}
 k(m,n)^{-1/2}
 \leq X^{-1/24}N^{1+o(1)}=o(N).
\]
The full GCD sum is at least $N$, since each diagonal pair $m=n$
contributes $1$.  The discarded tail is negligible
relative to that sum for every admissible family $\M$.
For every remaining pair, \eqref{eq:Wsmall}, with $\delta=1/4$, bounds
$W_{\ell,d}(k)$ below by a positive constant times $(\log X)^\ell$
throughout $X<d\leq2X$.
These square main terms give
\[
 S_{2,\ell}\gg_\ell
 Xh(Q)(\log X)^\ell
 \sum_{m,n\in\M}\sqrt{\frac{(m,n)}{[m,n]}}.
\]
The error \eqref{eq:central-error} is negligible here: after division
by $Xh(Q)N(\log X)^\ell$, it is at most
\[
 X^{-1/4+\eps+o(1)}N^{1+\eps}
 =X^{-\beta+\eps(5/4-\beta)+o(1)}=o(1)
\]
when $\eps$ is sufficiently small in terms of the fixed $\beta$.

We have $R_d^2\geq0$, and \eqref{eq:central-den} is positive for large
$X$.  The quotient $S_{2,\ell}/S_1$ is a weighted average of the real numbers
$(-1)^\ell L^{(\ell)}(1/2,\chi_d)$.
Their maximum is at least this average.
The lower bound for the kernel uses the fixed value $\delta=1/4$,
so its positive constant depends only on $\ell$.
Combining the numerator estimate with \eqref{eq:central-den}
gives, for some $c_\ell>0$,
\[
 \max_d(-1)^\ell L^{(\ell)}(\tfrac12,\chi_d)
 \geq (c_\ell+o(1))(\log X)^\ell
 \frac1N\sum_{m,n\in\mathcal M}
 \sqrt{\frac{(m,n)}{[m,n]}}.
\]
Both the discarded tail and the orthogonality error are
negligible relative to the displayed GCD sum.
The prime-support assumption was the only property of $\M$ used
in these error estimates.  The argument applies to every
family in the statement and proves \eqref{eq:central-precise}.
\end{proof}

\begin{proof}[Proof of Corollary~\ref{cor:central-exponent}]
For each fixed $0<\beta<1/4$, take
$N=\lfloor X^{1/4-\beta}\rfloor$ and choose the family $\M$
provided by \eqref{eq:GCD}--\eqref{eq:pmax}.
Its prime support satisfies the hypothesis of Theorem~\ref{thm1}.
Combining that theorem with \eqref{eq:GCD} gives
\[
 \max_d(-1)^\ell L^{(\ell)}(\tfrac12,\chi_d)
 \geq(c_\ell+o(1))(\log X)^\ell
 \exp\!\left((2+o(1))
 \sqrt{\frac{\log N\log_3N}{\log_2N}}\right).
\]
For fixed $\beta$,
$\log N=(1/4-\beta)\log X+o(1)$ and
$\log_jN=(1+o(1))\log_jX$ for $j=2,3$.
The exponential on the right has logarithm
$(2\sqrt{1/4-\beta}+o(1))\mathcal V_X$.
The other factors have logarithm
$\ell\log_2X+\log(c_\ell+o(1))=o(\mathcal V_X)$, because
$\log_2X=o(\mathcal V_X)$.
For each fixed $\beta$, this yields
\[
 \max_d(-1)^\ell L^{(\ell)}(\tfrac12,\chi_d)
 \geq\exp\bigl((2\sqrt{1/4-\beta}+o(1))\mathcal V_X\bigr).
\]
Let $\beta$ tend to zero through fixed positive values and choose
the corresponding thresholds for $X$ successively so that all
preceding estimates hold beyond each threshold.
This diagonal choice makes the coefficient tend to $1$ while
preserving the error bounds.  This proves the corollary.
\end{proof}

\section{The mesoscopic central band}\label{sec:mesoscopic}

We now prove Theorem~\ref{thm:mesoscopic}.  The first lemma describes
the kernels for a mesoscopic shift from the exact center.  We use
the functions $V_u$ and $X_u(d)$ defined in Section~\ref{sec3}.

\begin{lemma}\label{lem:shifted-kernel}
Let $0\leq\alpha<1/2$ and $X<d\leq2X$, and put
\begin{align*}
 W^+_{\ell,\alpha,d}(r)
  &=(-1)^\ell\left.\partial_u^\ell
    \{r^{-u}V_u(r/\sqrt d)\}\right|_{u=\alpha},\\
 W^-_{\ell,\alpha,d}(r)
  &=(-1)^\ell\left.\partial_u^\ell
    \{X_u(d)r^uV_{-u}(r/\sqrt d)\}\right|_{u=\alpha}.
\end{align*}
Then
\begin{equation}\label{eq:shifted-afe}
 (-1)^\ell L^{(\ell)}(\tfrac12+\alpha,\chi_d)
 =\sum_{r\geq1}\frac{\chi_d(r)}{\sqrt r}
 \{W^+_{\ell,\alpha,d}(r)+W^-_{\ell,\alpha,d}(r)\}.
\end{equation}
For every fixed $0<\eta<1/2$ there is $c_\eta>0$ such that, uniformly
for $r\leq X^\eta$,
\begin{equation}\label{eq:shifted-first-small}
 W^+_{\ell,\alpha,d}(r)
 =r^{-\alpha}(\log r)^\ell
 +O_{\ell,\eta}\!\left(
 X^{-c_\eta}r^{-\alpha}(\log X)^\ell\right).
\end{equation}
Moreover, for every $C>0$, one has
\begin{align}
 |W^+_{\ell,\alpha,d}(r)|
  &\ll_{\ell,C}r^{-\alpha}\log^\ell(2Xr)
       (1+r/\sqrt X)^{-C},\label{eq:shifted-first-global}\\
 |W^-_{\ell,\alpha,d}(r)|
  &\ll_{\ell,C}X^{-\alpha}r^\alpha\log^{\ell+1}(2Xr)
       (1+r/\sqrt X)^{-C}.\label{eq:shifted-dual-global}
\end{align}
The same estimates hold after applying $d\partial_d$.
For $0\leq\alpha\leq1/4$, the power $\ell+1$ of the logarithm in
\eqref{eq:shifted-dual-global} may be replaced by $\ell$,
also after applying $d\partial_d$.
\end{lemma}

\begin{proof}
The identity \eqref{eq:shifted-afe} is the differentiated form of
\cite[Lemma~2.2]{NearCenter}, with $G\equiv1$.  We verify below
the additional uniformity of the kernel estimates as $\alpha\uparrow1/2$.
Mellin inversion gives
\[
 V_u(x)=\frac{\Gamma(1/4+u/2,\pi x^2)}{\Gamma(1/4+u/2)}
 =\frac1{\Gamma(1/4+u/2)}
   \int_{\pi x^2}^{\infty}e^{-t}t^{-3/4+u/2}\dd t.
\]
For $0\leq u\leq1/2$, the incomplete-gamma formula implies
\[
 \partial_u^j\{V_u(x)-1\}
 \ll_j x^{1/2+u}(1+|\log x|)^j
 \qquad(0<x\leq1,\ 0\leq j\leq\ell).
\]
This follows by differentiating the integral over $0<t<\pi x^2$:
each parameter derivative inserts a power of $\tfrac12\log t$,
and the gamma denominator and its derivatives are uniformly bounded,
since $1/4+u/2\in[1/4,1/2]$.
For $x=r/\sqrt d\leq X^{-1/2+\eta}$, these bounds decay as a
fixed negative power of $X$, after absorbing logarithms into a
slightly smaller exponent.
Leibniz's rule applied to $r^{-u}V_u(x)$ proves
\eqref{eq:shifted-first-small}.

The same integral representation proves the global bound for the
first term.  Its differentiated integrand is bounded by
\[
 r^{-\alpha}e^{-t}t^{-3/4+\alpha/2}
 O_\ell\bigl((1+\log r+|\log t|)^\ell\bigr).
\]
For a lower limit at most $1$, its integral is uniformly bounded
by $O_\ell(r^{-\alpha}\log^\ell(2Xr))$.
For a lower limit at least $1$, exponential decay in $t$ gives
any prescribed negative power of $r/\sqrt X$.
This proves \eqref{eq:shifted-first-global}.

For the dual term, cancel the factor $\Gamma(1/4-u/2)$ between
$X_u(d)$ and $V_{-u}$ before differentiating.  The result is
\[
 \frac1{\Gamma(1/4+u/2)}
 \int_{\pi r^2/d}^{\infty}
 e^{-t}t^{-3/4}
 \exp\!\left(-u\left\{\log\frac d{\pi r}
                    +\frac12\log t\right\}\right)\dd t.
\]
There is no singular gamma factor as $u$ tends to $1/2$.
After factoring out $(d/\pi)^{-u}r^u$, the remaining power of $t$
is $t^{-3/4-u/2}$.
If $z=\pi r^2/d\leq1$, then uniformly for $0\leq u\leq1/2$,
\[
 \int_z^1 t^{-3/4-u/2}(1+|\log t|)^\ell\dd t
 \leq\int_z^1t^{-1}(1+|\log t|)^\ell\dd t
 \ll_\ell(1+|\log z|)^{\ell+1}.
\]
Here $|\log z|\ll\log(2Xr)$.
The part with $t\geq1$ is bounded by exponential decay.
In Leibniz's rule, a term with $j$ derivatives
on the integral has at most $\ell-j$ logarithms from the prefactor
and $j+1$ from the integral.  Its total degree is at most $\ell+1$.
If $z>1$, that decay instead gives an
arbitrary negative power of $r/\sqrt X$.
Since $(d/\pi)^{-\alpha}\asymp X^{-\alpha}$ uniformly in the
stated range, this proves \eqref{eq:shifted-dual-global}.
Finally, applying $d\partial_d$ differentiates the elementary
prefactor and the lower limit.  The latter contributes a factor
of $z$ times the integrand at $z$. For the dual integral this
contains $e^{-z}z^{1/4-u/2}$, which is bounded when $z\leq1$
and decays exponentially when $z>1$.
The same argument applies to the first integral, and proves the
claimed derivative estimates.

For the sharper dual estimate, $0\leq u\leq1/4$ gives
$1/4-u/2\geq1/8$, and consequently
\[
 \int_0^1t^{-3/4-u/2}(1+|\log t|)^\ell\dd t\ll_\ell1
\]
uniformly in $u$.  The integral over $z\leq t\leq1$
does not require the extra power of $\log(2Xr)$.
The endpoint terms after $d\partial_d$ have
$z^{1/4-u/2}(1+|\log z|)^\ell\ll_\ell1$ for $0<z\leq1$.
The range $z>1$ still has exponential decay.
This proves the sharper bound, including its $d$-derivative form.
\end{proof}

\begin{lemma}\label{lem:shifted-positive}
Fix $\ell\geq0$.  There is an integer $R_\ell\geq1$ such that, uniformly
for $0\leq\alpha<1/2$ and all sufficiently large $X$,
\begin{align}
 W^-_{\ell,\alpha,d}(r)&\geq0\qquad(r\geq1),\label{eq:dual-positive}\\
 W^+_{\ell,\alpha,d}(r)&\geq0\qquad(r\geq R_\ell),\label{eq:first-positive}\\
 W^+_{\ell,\alpha,d}(r)&\geq-O_\ell(1)\qquad(1\leq r<R_\ell).
 \label{eq:first-finite-negative}
\end{align}
If $\varepsilon_X\to0$ and $r\leq X^{\varepsilon_X}$, then, uniformly
when $\alpha=o(1)$,
\begin{align}
 W^+_{\ell,\alpha,d}(r)
 &=r^{-\alpha}(\log r)^\ell+
 o\bigl(r^{-\alpha}(1+\log r)^\ell\bigr),\label{eq:first-boundary}\\
 W^-_{\ell,\alpha,d}(r)
 &=(1+o(1))X^{-\alpha}r^\alpha(\log X)^\ell.
 \label{eq:dual-boundary}
\end{align}
In \eqref{eq:first-boundary}, the assertion is used only when
$\log r\to\infty$ if $\ell>0$.
\end{lemma}

\begin{proof}
For $0\leq u\leq1/2$ and real $q$, define
\[
 \mathfrak P_{\ell,u}(q):=
 (-1)^\ell\Gamma(1/4+u/2)\mathrm e^{uq}
 \partial_u^\ell
 \left\{\frac{\mathrm e^{-uq}}{\Gamma(1/4+u/2)}\right\}.
\]
This is a polynomial in $q$.  Its degree and uniform coefficient
bounds can be checked directly from
\[
 \mathfrak P_{0,u}(q)=1,\qquad
 \mathfrak P_{\ell+1,u}(q)
 =\bigl(q+\tfrac12\psi(1/4+u/2)\bigr)
       \mathfrak P_{\ell,u}(q)-\partial_u\mathfrak P_{\ell,u}(q).
\]
The recurrence follows by differentiating
$\mathrm e^{-uq}/\Gamma(1/4+u/2)$ once more.
All fixed-order polygamma functions are bounded on $[1/4,1/2]$,
so induction shows that the leading coefficient is $1$ and all
other coefficients are bounded uniformly in $u$.
For $\ell\geq1$ this proves
\begin{equation}\label{eq:Bell-q}
 \mathfrak P_{\ell,u}(q)=q^\ell+O_\ell((1+q)^{\ell-1})
 \quad(q\geq0),\qquad
 |\mathfrak P_{\ell,u}(q)|\ll_\ell(1+|q|)^\ell.
\end{equation}
For $\ell=0$ the polynomial is exactly $1$.
Choose a fixed $Q_\ell\geq1$ so that
$\mathfrak P_{\ell,u}(q)\geq\tfrac12q^\ell$ whenever
$q\geq Q_\ell$, uniformly in $u$.

Put $z=\pi r^2/d$ and
\[
 q_+(t,r)=\log r-\tfrac12\log t,\qquad
 q_-(t,r,d)=\log\frac d{\pi r}+\tfrac12\log t.
\]
Differentiating the two integral formulas above gives the exact
identities
\begin{align}
 W^\pm_{\ell,\alpha,d}(r)
 =\frac1{\Gamma(1/4+\alpha/2)}
 \int_z^\infty \mathrm e^{-t}t^{-3/4}
 \mathrm e^{-\alpha q_\pm(t)}
 \mathfrak P_{\ell,\alpha}(q_\pm(t))\,\dd t.
 \label{eq:Laplace-two-halves}
\end{align}
Since $t\geq z$,
\[
 q_-(t,r,d)\geq
 \log\frac d{\pi r}+\frac12\log\frac{\pi r^2}d
 =\frac12\log\frac d\pi.
\]
For large $X$ this is greater than $Q_\ell$, independently of $r$.
Every factor in the dual integrand is then nonnegative, proving
\eqref{eq:dual-positive}.

For the first half, $q_+(t,r)\geq Q_\ell$ when
$t\leq T_r:=r^2\mathrm e^{-2Q_\ell}$.
Any negative contribution must come from $t>T_r$.
For $T_r\geq2$ and $t>T_r$, one has
$\log r\leq Q_\ell+\tfrac12\log t$.
Using \eqref{eq:Bell-q} and
$\mathrm e^{-\alpha q_+}=r^{-\alpha}t^{\alpha/2}$,
the absolute value of this tail is bounded uniformly in $\alpha$ by
\[
 O_\ell\!\left(\int_{T_r}^\infty
 e^{-t}t^{-1/2}(1+\log t)^\ell\dd t\right)
 =O_\ell(e^{-T_r/2}).
\]
If $z\leq1$, the interval $1\leq t\leq2$ lies inside the integration
range.  For sufficiently large fixed $r$, it also lies below $T_r$,
and its positive contribution is
\[
 \gg_\ell r^{-\alpha}(\log r)^\ell
 \geq r^{-1/2}(\log r)^\ell.
\]
Since $T_r$ is a fixed positive multiple of $r^2$,
$e^{-T_r/2}=o(r^{-1/2}(\log r)^\ell)$ as $r$ tends to infinity.
The positive integral dominates the possibly negative
tail for every $r$ beyond a fixed threshold.

If $z>1$, use the interval $z\leq t\leq z+1$.
For large $d$ this lies below $T_r$, because
\[
 \frac{T_r}{z}=\frac{d\mathrm e^{-2Q_\ell}}\pi\longrightarrow\infty.
\]
On this interval,
$q_+(t,r)\geq\tfrac12\log(d/(2\pi))$.  We also have
$r^{-\alpha}t^{\alpha/2}\geq(\pi/d)^{\alpha/2}
\gg d^{-1/4}$.
Its positive contribution is
\[
 \gg_\ell e^{-z}z^{-3/4}d^{-1/4}(\log d)^\ell.
\]
The ratio of the preceding tail bound to this lower bound is at most
\[
 O_\ell\!\left(
 d^{1/4}z^{3/4}(\log d)^{-\ell}
 \exp\!\left\{-\left(
 \frac{d\mathrm e^{-2Q_\ell}}{2\pi}-1\right)z\right\}\right)=o(1),
\]
uniformly for $z>1$.  The exponential dominates both powers, since its decay rate tends
to infinity with $d$.
Choose the fixed threshold $R_\ell$ to cover the first case.
Taking $X$ sufficiently large then covers both cases and proves
\eqref{eq:first-positive}.
For the finitely many integers $1\leq r<R_\ell$,
\eqref{eq:shifted-first-small} shows that the first weight equals
the nonnegative quantity $r^{-\alpha}(\log r)^\ell$ plus a bounded
error, uniformly in $\alpha$.  This proves
\eqref{eq:first-finite-negative}.

Finally, if $r\leq X^{\varepsilon_X}$, then
$z=X^{-1+o(1)}$.
For $\alpha=o(1)$, the gamma parameters $1/4\pm\alpha/2$ remain
in a fixed compact subinterval of $(0,\infty)$.
The lower incomplete-gamma integrals and all their fixed parameter
derivatives are therefore $X^{-c+o(1)}$ for some $c>0$.
Leibniz's rule for the first prefactor gives
\eqref{eq:first-boundary}.  Its error is smaller than
$r^{-\alpha}(1+\log r)^\ell$ even after multiplication by fixed
powers of $\log X$.

For the dual prefactor, logarithmic differentiation gives
\[
 -\partial_u\log\{X_u(d)r^u\}
 =\log\frac d{\pi r}
 +\frac12\psi(1/4-u/2)+\frac12\psi(1/4+u/2)
 =(1+o(1))\log X
\]
at $u=\alpha$, uniformly for the indicated $d$ and $r$.
Here $\log r=o(\log X)$ and both polygamma arguments stay away
from zero.  All higher logarithmic derivatives are $O_\ell(1)$.
The leading term in the $\ell$th derivative is
$(1+o(1))(\log X)^\ell$ times the prefactor.  Every term
involving a higher logarithmic derivative is of smaller order.
Also
\[
 X_\alpha(d)r^\alpha
 =X^{-\alpha}r^\alpha
   \exp\!\left\{-\alpha\log\frac d{\pi X}
       +\log\frac{\Gamma(1/4-\alpha/2)}
                      {\Gamma(1/4+\alpha/2)}\right\}
 =(1+o(1))X^{-\alpha}r^\alpha.
\]
The derivatives of $V_{-u}-1$ contribute only a relative
$X^{-c+o(1)}$ error times fixed powers of $\log X$.
This proves \eqref{eq:dual-boundary}.
\end{proof}

The new ingredient is the following near-critical, log-weighted
G\'al-sum construction.  Its normalization is chosen so that the
quadratic mass in each prime block does not lose an exponential factor
when $\alpha\log_2N$ is nonzero.  The block decomposition originates in
the long-resonator constructions of Bondarenko--Seip and Chirre
\cite{BondarenkoSeip,Chirre}.  For its recent quadratic
adaptation, see \cite{NearCenter}.

\begin{lemma}
\label{lem:mesoscopic-resonator}
Fix $A_0>0$, allow $0\leq A\leq A_0$ to depend on $N$, and fix
$0<\gamma<1$, $a>1$ with $a\gamma<1$, and $\ell\geq0$.
Let $N\to\infty$, and put
\[
 L=\log N,\quad L_2=\log_2N,\quad L_3=\log_3N,\quad
 Y_0=LL_2,\quad K=\lfloor L_2^\gamma\rfloor,
\]
and suppose $\alpha\geq0$ and
\[
 \alpha=\frac{A+o(1)}{L_2}.
\]
For $1\leq j\leq K$ put
\[
 \mathcal P_j=\{p:\ \mathrm e^jY_0<p\leq \mathrm e^{j+1}Y_0\},
\]
and define the square-free multiplicative function $f$ on
$\mathcal P:=\bigcup_{j\leq K}\mathcal P_j$ by
\begin{equation}\label{eq:mesoscopic-f}
 \begin{aligned}
 f(p)&=\frac{C_N}{p^{1/2+\alpha}\log(p/Y_0)}
       \qquad(p\in\mathcal P),\\
 C_N&=Y_0^{2\alpha}
 \frac{L^{1/2-\alpha}L_2^{1/2+\alpha}}{L_3^{1/2-\alpha}},
 \end{aligned}
\end{equation}
with $f(p)=0$ otherwise.  Let $\mathcal M$ consist of those
$m\in\operatorname {supp}(f)$ for which
\begin{equation}\label{eq:mesoscopic-count-cutoff}
 \omega_{\mathcal P_j}(m)\leq
 \Delta_j:=\left\lfloor\frac{aL}{L_3j(j+1)}\right\rfloor
 \qquad(1\leq j\leq K).
\end{equation}
Then $|\mathcal M|\leq N$.  Put
\[
 k(m,n)=\frac{mn}{(m,n)^2},\qquad
 H_N=\sum_{m\in\mathcal M}f(m)^2h(m),
\]
and define
\begin{align*}
 G_N&=\sum_{p\in\mathcal P}
 \log\frac{1+h(p)f(p)^2+2h(p)f(p)p^{-1/2-\alpha}}
 {1+h(p)f(p)^2},\\
 \mu_N&=\sum_{p\in\mathcal P}
 \frac{2h(p)f(p)p^{-1/2-\alpha}}
 {1+h(p)f(p)^2+2h(p)f(p)p^{-1/2-\alpha}}\log p.
\end{align*}
Also put
\begin{align*}
 G_N^-&=\sum_{p\in\mathcal P}
 \log\frac{1+h(p)f(p)^2+2h(p)f(p)p^{-1/2+\alpha}}
 {1+h(p)f(p)^2}.
\end{align*}
Then, for every fixed $\varepsilon>0$,
\begin{align}
 &\sum_{\substack{m,n\in\mathcal M\\k(m,n)\leq N^\varepsilon}}
 \frac{f(m)f(n)h([m,n])}{k(m,n)^{1/2+\alpha}}
 \bigl(\log k(m,n)\bigr)^\ell \notag\\
 &\qquad=(1+o(1))H_N\mu_N^\ell\exp(G_N),
 \label{eq:mesoscopic-gal}
\end{align}
where the factor $\mu_N^\ell$ is omitted when $\ell=0$, and
\begin{align}
 G_N&=(2\gamma\mathrm e^{-A}+o(1))
       \sqrt{\frac{L L_3}{L_2}},\label{eq:mesoscopic-free-energy}\\
 \mu_N&=(2\gamma\mathrm e^{-A}+o(1))
       \sqrt{L L_2L_3}.
 \label{eq:mesoscopic-mean}
\end{align}
If $\varepsilon_N=(L_2L_3/L)^{1/4}$, then
\begin{align}
 &\sum_{\substack{m,n\in\mathcal M\\
 k(m,n)\leq N^{\varepsilon_N}}}
 \frac{f(m)f(n)h([m,n])}{k(m,n)^{1/2-\alpha}}
 =(1+o(1))H_N\exp(G_N^-),\label{eq:opposite-general}\\
 &G_N^-=(2\gamma\mathrm e^A+o(1))
 \sqrt{\frac{LL_3}{L_2}}.\label{eq:opposite-general-energy}
\end{align}
All displayed estimates are uniform for $0\leq A\leq A_0$.
\end{lemma}

\begin{proof}
We first compute the three prime sums that govern the construction.
For $1\leq t\leq K+1$ we have
\[
 \alpha t=o(1),\qquad
 \log(Y_0\mathrm e^t)=L_2+L_3+t=L_2(1+o(1)),
\]
uniformly, because $K=L_2^{\gamma}+O(1)$ and $\gamma<1$.
The prime number theorem and partial summation give
\begin{align*}
 \sum_{p\in\mathcal P_j}f(p)^2
 &=(1+o(1))C_N^2
   \int_{Y_0\mathrm e^j}^{Y_0\mathrm e^{j+1}}
   \frac{\dd x}{x^{1+2\alpha}\log^2(x/Y_0)\log x}\\
 &=(1+o(1))\frac{C_N^2Y_0^{-2\alpha}}{L_2}
   \int_j^{j+1}\frac{\mathrm e^{-2\alpha t}}{t^2}\dd t,\\
 \sum_{p\in\mathcal P_j}\frac{f(p)}{p^{1/2+\alpha}}
 &=(1+o(1))\frac{C_NY_0^{-2\alpha}}{L_2}
   \int_j^{j+1}\frac{\mathrm e^{-2\alpha t}}t\dd t.
\end{align*}
The prime-number-theorem error is uniform over all these intervals.
Their endpoints tend to infinity, each interval has a fixed ratio of
endpoints, and the logarithmic derivatives of the displayed weights
are bounded on each interval.  The normalization gives
\begin{align*}
 \frac{C_N^2Y_0^{-2\alpha}}{L_2}
 &=\frac L{L_3}L_2^{4\alpha}L_3^{2\alpha}
   =(1+o(1))\frac L{L_3},\\
 \frac{C_NY_0^{-2\alpha}}{L_2}
 &=\sqrt{\frac{L}{L_2L_3}}
   \exp\{-\alpha L_2+\alpha L_3+\alpha\log L_3\}
   =(\mathrm e^{-A}+o(1))\sqrt{\frac{L}{L_2L_3}}.
\end{align*}
Multiplying the second summand by $\log p=(1+o(1))L_2$
gives the third prime sum.  We conclude, uniformly for $1\leq j\leq K$,
\begin{align}
 \sum_{p\in\mathcal P_j}f(p)^2
  &=(1+o(1))\frac{L}{L_3}
    \int_j^{j+1}\frac{\dd t}{t^2}
   =(1+o(1))\frac{L}{L_3j(j+1)},
 \label{eq:block-square-mass}\\
 \sum_{p\in\mathcal P_j}\frac{f(p)}{p^{1/2+\alpha}}
  &=(\mathrm e^{-A}+o(1))
    \sqrt{\frac{L}{L_2L_3}}
    \int_j^{j+1}\frac{\dd t}{t},
 \label{eq:block-linear-mass}\\
 \sum_{p\in\mathcal P_j}\frac{f(p)\log p}{p^{1/2+\alpha}}
  &=(\mathrm e^{-A}+o(1))
    \sqrt{\frac{LL_2}{L_3}}
    \int_j^{j+1}\frac{\dd t}{t}.
 \label{eq:block-log-mass}
\end{align}
The factor $Y_0^{2\alpha}$ in \eqref{eq:mesoscopic-f} is essential:
omitting it multiplies the square mass by $Y_0^{-4\alpha}
=\mathrm e^{-4A+o(1)}$ and changes the occupancy budget.  Summing
\eqref{eq:block-linear-mass} and \eqref{eq:block-log-mass}, and using
\[
 \sum_{j\leq K}\int_j^{j+1}\frac{\dd t}{t}
 =\log(K+1)=\gamma L_3+o(L_3),
\]
we obtain
\begin{align}
 \sum_{p\in\mathcal P}\frac{f(p)}{p^{1/2+\alpha}}
  &=(\gamma\mathrm e^{-A}+o(1))
    \sqrt{\frac{LL_3}{L_2}},\label{eq:linear-total}\\
 \sum_{p\in\mathcal P}\frac{f(p)\log p}{p^{1/2+\alpha}}
  &=(\gamma\mathrm e^{-A}+o(1))
    \sqrt{LL_2L_3}.
 \label{eq:log-total}
\end{align}

We next verify that the block cutoffs give at most $N$ integers.
An integer in $\mathcal M$ is determined by choosing at most $\Delta_j$
primes from each $\mathcal P_j$.  The prime number theorem gives
$|\mathcal P_j|\leq2\mathrm e^{j+1}L$ for large $N$.  It also implies
$\Delta_j<|\mathcal P_j|/2$ uniformly in $j$.  The binomial coefficients
are increasing up to this cutoff, so
\[
 |\mathcal M|
 \leq\prod_{j\leq K}\sum_{v=0}^{\Delta_j}
       \binom{|\mathcal P_j|}{v}
 \leq\prod_{j\leq K}(\Delta_j+1)
       \left(\frac{\mathrm e|\mathcal P_j|}{\Delta_j}\right)^{\Delta_j}.
\]
The minimum of $L/(L_3j(j+1))$ tends to infinity.  Replacing this
quantity by its floor changes $\Delta_j$ by a relative $o(1)$ uniformly,
and
\[
 \log\frac{\mathrm e|\mathcal P_j|}{\Delta_j}
 \leq j+2\log(j+1)+\log L_3+O_a(1).
\]
The factors $\Delta_j+1$ contribute at most $O(K\log L)=o(L)$
to the logarithm.  Therefore
\begin{align*}
 \log|\mathcal M|
 &\leq\frac{aL}{L_3}\sum_{j\leq K}
  \frac{j+2\log(j+1)+\log L_3+O_a(1)}{j(j+1)}+o(L)\\
 &=\frac{aL}{L_3}\{\log K+O_a(\log L_3+1)\}+o(L)
  =(a\gamma+o(1))L<L.
\end{align*}
In the last line we used $\sum_{j\leq K}1/(j+1)=\log K+O(1)$,
$\sum_{j\geq1}\log(j+1)/(j(j+1))<\infty$, and
$\log K=\gamma L_3+o(1)$.  Hence $|\mathcal M|\leq N$.

The same cutoff gives
\[
 \omega(m)\leq\sum_{j\leq K}\Delta_j\leq\frac{aL}{L_3},
 \qquad
 \log m\leq\sum_{j\leq K}\Delta_j(\log Y_0+j+1)
 \ll_{a,\gamma}\frac{LL_2}{L_3}.
\]
We record the latter estimate as
\begin{equation}\label{eq:mesoscopic-element-size}
 \log m\ll_{a,\gamma}\frac{LL_2}{L_3}\qquad(m\in\mathcal M).
\end{equation}
Every prime divisor of $m$ exceeds $\mathrm eY_0$, so
\[
 0\leq-\log h(m)
 \leq\sum_{p\mid m}\frac1p
 \ll_a\frac1{L_2L_3}=o(1).
\]
We have $h(m)=1+o(1)$ uniformly on $\mathcal M$.  We use
\eqref{eq:mesoscopic-element-size} to bound the orthogonality errors
below.

It remains to evaluate the quadratic form.  First sum over the full
support of $f$, and introduce the moment generating product
\begin{equation}\label{eq:mesoscopic-Z}
 Z_N(z)=\prod_{p\in\mathcal P}
 \left(1+h(p)f(p)^2+2h(p)f(p)p^{-1/2-\alpha+z}\right).
\end{equation}
The three local states are: $p$ divides both $m,n$, $p$ divides exactly
one of them (with two possible orientations), or $p$ divides neither.
The derivative $Z_N^{(\ell)}(0)$ is the full version of the sum on the
left of \eqref{eq:mesoscopic-gal}.  Put
\[
 H_N^*=\prod_{p\in\mathcal P}(1+h(p)f(p)^2).
\]
The square-mass calculation also gives
$\max_pf(p)^2\ll L_3^{-1}$, whereas $h(p)=1+O(p^{-1})$.
In each local quotient defining $G_N$, the denominator is $1+o(1)$
uniformly.  Its added linear term also tends to zero uniformly.  The
expansion $\log(1+x)=x+O(x^2)$ gives
\begin{equation}\label{eq:Z-free-energy}
 \log\frac{Z_N(0)}{H_N^*}
 =G_N=(2+o(1))\sum_{p\in\mathcal P}\frac{f(p)}{p^{1/2+\alpha}}.
\end{equation}
For example, the total quadratic error is bounded by the largest
linear term times their sum, so it is $o(1)$ times that sum.
This justifies the relative error in \eqref{eq:Z-free-energy}.

After division by $Z_N(0)$, the quantity $S=\log k(m,n)$ is a sum
of independent variables $\xi_p\log p$, where $\xi_p$ is Bernoulli with
parameter
\[
 q_p=\frac{2h(p)f(p)p^{-1/2-\alpha}}
 {1+h(p)f(p)^2+2h(p)f(p)p^{-1/2-\alpha}}.
\]
We have $\mathbb E S=\mu_N$, and $\log p\ll L_2$ throughout
$\mathcal P$.  For each fixed integer $j\geq2$, the $j$-th cumulant
of $S$ satisfies
\[
 |\kappa_j(S)|\ll_j\sum_{p\in\mathcal P}q_p(\log p)^j
 \ll_j L_2^{j-1}\mu_N.
\]
Since $\mu_N\asymp\sqrt{LL_2L_3}$ uniformly for bounded $A$,
\[
 \frac{|\kappa_j(S)|}{\mu_N^j}
 \ll_j\left(\frac{L_2}{\mu_N}\right)^{j-1}=o(1)
 \qquad(j\geq2\text{ fixed}).
\]
The moment--cumulant argument in the proof of
Lemma~\ref{lem:free-energy} gives
$\mathbb E S^j=(1+o(1))\mu_N^j$ for every fixed integer $j\geq0$,
uniformly in the present parameters. The zeroth moment is $1$.
In particular,
\begin{equation}\label{eq:Z-log-moment}
 Z_N^{(\ell)}(0)=(1+o(1))\mu_N^\ell Z_N(0).
\end{equation}

We now show that the cutoffs preserve this relative asymptotic, including
the logarithmic weight.  Write $b_j=L/(L_3j(j+1))$.  Under the probability
measure defined by \eqref{eq:mesoscopic-Z}, the number of primes of
$\mathcal P_j$ dividing either fixed member of the pair is a sum of
independent Bernoulli variables with mean
\[
 \sum_{p\in\mathcal P_j}
 \frac{h(p)f(p)^2+h(p)f(p)p^{-1/2-\alpha}}
 {1+h(p)f(p)^2+2h(p)f(p)p^{-1/2-\alpha}}
 =(1+o(1))b_j.
\]
Here the one-sided contribution is $o(b_j)$ uniformly in $j$, since
\[
 \frac{p^{-1/2-\alpha}}{f(p)}
 =\frac{\log(p/Y_0)}{C_N}\ll\frac K{C_N}=o(1).
\]
For a sum of independent Bernoulli variables with mean
$(1+o(1))b_j$, exponential Markov with parameter $\log a>0$ gives
\begin{align*}
 \mathbb P\{\omega_{\mathcal P_j}(m)>\Delta_j\}
 &\leq\exp\{-(\Delta_j+1)\log a+(a-1)(1+o(1))b_j\}\\
 &\leq\exp\{-\tfrac12(a\log a-a+1)b_j\}.
\end{align*}
The Bernoulli moment-generating function is at most
$\exp((a-1)\mathbb E\omega_{\mathcal P_j}(m))$, and
$\Delta_j+1\geq ab_j$.  The exponent has a strictly negative constant
because $a\log a-a+1>0$ for $a>1$.
The same estimate holds for $n$.  If $\mathcal B_{\rm occ}$ is the
event that either member lies outside $\mathcal M$, the union bound yields
\[
 \mathbb P(\mathcal B_{\rm occ})
 \leq2\sum_{j\leq K}
  \exp\{-\tfrac12(a\log a-a+1)b_j\}=o(1).
\]
Here $\min_{j\leq K}b_j\asymp L/(L_3L_2^{2\gamma})$ grows faster
than $\log K$, which makes the whole sum tend to zero.
Under the product measure defining $H_N^*$, the corresponding Bernoulli
parameters are $h(p)f(p)^2/(1+h(p)f(p)^2)$. Their block means also equal
$(1+o(1))b_j$.  The same tail estimate yields
\[
 H_N=(1+o(1))H_N^*.
\]

Let $\eta_N$ be either a fixed positive $\varepsilon$ or
$\varepsilon_N=(L_2L_3/L)^{1/4}$.  Independence and $\log p/L_2=O(1)$ give
\[
 \log\frac{Z_N(1/L_2)}{Z_N(0)}
 =\sum_{p\in\mathcal P}\log\bigl(1+q_p(e^{\log p/L_2}-1)\bigr)
 \ll\frac1{L_2}\sum_{p\in\mathcal P}q_p\log p
 =O(\mu_N/L_2).
\]
For $\mathcal B_k=\{S>\eta_NL\}$,
\[
 \mathbb P(\mathcal B_k)
 \leq\exp\left(-\frac{\eta_NL}{L_2}
                  +O(\mu_N/L_2)\right)=o(1).
\]
For the shrinking cutoff this follows from
\[
 \frac{\varepsilon_NL}{\mu_N}
 \asymp\left(\frac{L}{L_2L_3}\right)^{1/4}\longrightarrow\infty;
\]
the assertion for fixed $\varepsilon$ follows a fortiori.
Put $\mathcal B=\mathcal B_{\rm occ}\cup\mathcal B_k$.
By Cauchy--Schwarz and the preceding moment estimate with $j=2\ell$,
\[
 \mathbb E\bigl[S^\ell\mathbf1_{\mathcal B}\bigr]
 \leq(\mathbb E S^{2\ell})^{1/2}\mathbb P(\mathcal B)^{1/2}
 =o(\mu_N^\ell).
\]
This includes $\ell=0$ under the zeroth-moment convention.
Deleting the bad pairs changes $Z_N^{(\ell)}(0)$ by
$o(\mu_N^\ell Z_N(0))$.  Together with $H_N=(1+o(1))H_N^*$ and
\eqref{eq:Z-free-energy}, this proves
\eqref{eq:mesoscopic-gal}--\eqref{eq:mesoscopic-mean}.  The same argument
proves \eqref{eq:mesoscopic-gal} with $\varepsilon_N$ in place of fixed
$\varepsilon$.

For the opposite tilt, replace $p^{-1/2-\alpha}$ in
\eqref{eq:mesoscopic-Z} by $p^{-1/2+\alpha}$.  In this case
$f(p)p^{-1/2+\alpha}=C_N/(p\log(p/Y_0))$, so the corresponding
block sum is
\[
 \sum_{p\in\mathcal P_j}f(p)p^{-1/2+\alpha}
 =(1+o(1))\frac{C_N}{L_2}\int_j^{j+1}\frac{\dd t}{t}.
\]
The coefficient equals
\[
 \frac{C_N}{L_2}
 =\sqrt{\frac{L}{L_2L_3}}
  \exp\{\alpha L_2+3\alpha L_3+\alpha\log L_3\}
 =(\mathrm e^A+o(1))\sqrt{\frac{L}{L_2L_3}}.
\]
Summing over the blocks and expanding the local logarithms as above
proves \eqref{eq:opposite-general-energy}.
The associated mean of $S$ is
\[
 \sum_{p\in\mathcal P}
 \frac{2h(p)f(p)p^{-1/2+\alpha}}
 {1+h(p)f(p)^2+2h(p)f(p)p^{-1/2+\alpha}}\log p
 =(2\gamma\mathrm e^A+o(1))\sqrt{LL_2L_3}.
\]
The one-sided block contribution is still $o(b_j)$, because
\[
 \frac{p^{-1/2+\alpha}}{f(p)}
 =\frac{p^{2\alpha}\log(p/Y_0)}{C_N}=o(1)
 \qquad(p\in\mathcal P)
\]
uniformly on $\mathcal P$.  The same block-tail and $k$-tail estimates
apply with this mean in place of $\mu_N$.  Taking the zeroth moment
proves \eqref{eq:opposite-general} with the shrinking cutoff.
\end{proof}

\begin{corollary}\label{cor:opposite-tilts}
In the notation of Lemma~\ref{lem:mesoscopic-resonator}, suppose in
addition that
\[
 \alpha\sqrt{LL_2L_3}\longrightarrow0,
\qquad
 \varepsilon_N:=\left(\frac{L_2L_3}{L}\right)^{1/4}.
\]
The condition is imposed on the actual shift $\alpha$, including the
error allowed in $\alpha L_2=A+o(1)$.  A condition only on $A$ would
not control the relative error of the exponential factors.
Put
\[
 G_N^0=\sum_{p\in\mathcal P}
 \log\frac{1+h(p)f(p)^2+2h(p)f(p)p^{-1/2}}
 {1+h(p)f(p)^2}.
\]
Then
\begin{equation}\label{eq:central-free-energy}
 G_N^0=(2\gamma+o(1))\sqrt{\frac{LL_3}{L_2}},
\end{equation}
and
\begin{align}
 &\sum_{\substack{m,n\in\mathcal M\\
 k(m,n)\leq N^{\varepsilon_N}}}
 \frac{f(m)f(n)h([m,n])}{k(m,n)^{1/2+\alpha}}
 (\log k(m,n))^\ell \notag\\
 &\qquad=(1+o(1))H_N
 (2\gamma\sqrt{LL_2L_3})^\ell\exp(G_N^0),
 \label{eq:boundary-positive-tilt}\\
 &\sum_{\substack{m,n\in\mathcal M\\
 k(m,n)\leq N^{\varepsilon_N}}}
 \frac{f(m)f(n)h([m,n])}{k(m,n)^{1/2-\alpha}}
 =(1+o(1))H_N\exp(G_N^0).
 \label{eq:boundary-negative-tilt}
\end{align}
\end{corollary}

\begin{proof}
We compare the free energies $G_N$ and $G_N^-$ from the lemma
with $G_N^0$.  The additional hypothesis implies
$\alpha L_2\to0$, and hence $A=o(1)$ by the relation
$\alpha L_2=A+o(1)$.  By the mean-value theorem applied to the local
factors,
\[
 |G_N-G_N^0|+|G_N^--G_N^0|
 \ll\alpha\sum_{p\in\mathcal P}
      \frac{f(p)\log p}{\sqrt p}\,p^\alpha
 \ll\alpha\sqrt{LL_2L_3}
 =\alpha L_2\sqrt{\frac{LL_3}{L_2}}=o(1).
\]
Here $p^\alpha=1+o(1)$ uniformly on $\mathcal P$, since
$\log p=(1+o(1))L_2$.  The prime-sum estimates in the lemma, now with
$A=o(1)$, yield \eqref{eq:central-free-energy} and
$\mu_N=(2\gamma+o(1))\sqrt{LL_2L_3}$.
The shrinking-cutoff conclusions of
Lemma~\ref{lem:mesoscopic-resonator} apply.  Combining them with
$\exp(G_N)=\exp(G_N^-)=(1+o(1))\exp(G_N^0)$ and the displayed
asymptotic for $\mu_N$ gives
\eqref{eq:boundary-positive-tilt} and
\eqref{eq:boundary-negative-tilt}.
\end{proof}

\begin{proof}[Proof of Theorem~\ref{thm:mesoscopic}]
Fix $0<\beta<1/4$, $0<\gamma<1$, and $a>1$ with
$a\gamma<1$, and take
\[
 N=\lfloor X^{1/4-\beta}\rfloor,\qquad
 \alpha=\frac A{\log_2X}.
\]
Construct $f$ and $\mathcal M$ by Lemma~\ref{lem:mesoscopic-resonator},
and set
\[
 R_d=\sum_{m\in\mathcal M}f(m)\chi_d(m).
\]
Expanding $R_d^2$ gives a sum of $f(m)f(n)\chi_d(mn)$.
Since $m,n$ are square-free, $mn$ is a square precisely when $m=n$.
Lemma~\ref{lem:orth} gives the diagonal main term
$2X\sum_{m\in\mathcal M}f(m)^2h(m)/\pi^2=2XH_N/\pi^2$.
We claim that the total orthogonality error is $o(XH_N)$.  This gives
\begin{equation}\label{eq:mesoscopic-den}
 \sum_{\substack{X<d\leq2X\\d\in\Fplus}}R_d^2
 =\left(\frac2{\pi^2}+o(1)\right)XH_N.
\end{equation}
To verify the claim, fix $0<\delta<\beta$ in Lemma~\ref{lem:orth}.
For $mn=n_0n_1^2$, \eqref{eq:mesoscopic-element-size} implies
\[
 \log g_1(n_0)\ll_\delta
       (LL_2/L_3)^{1-\delta}=o(\log X).
\]
The primes of $n_1$ belong to $\mathcal P$, so
\[
 \log g_2(n_1)
 \leq\sum_{p\mid mn}p^{-1/2-\delta}
 \leq\frac{2aL/L_3}{(\mathrm eY_0)^{1/2+\delta}}
 =o(\log X).
\]
Both arithmetic factors are $X^{o(1)}$ uniformly.  Cauchy's
inequality and $|\mathcal M|\leq N$ give
$(\sum_{m\in\mathcal M}f(m))^2\leq N\sum_{m\in\mathcal M}f(m)^2$.
Since $h(m)=1+o(1)$ uniformly, the total error is at most
\[
 X^{1/2+\delta+o(1)}N\sum_{m\in\mathcal M}f(m)^2
 =(1+o(1))X^{1/2+\delta+o(1)}NH_N=o(XH_N).
\]
The last step uses $N\leq X^{1/4-\beta}$. Its exponent relative to
$XH_N$ is $-1/4-\beta+\delta+o(1)<0$.

Insert \eqref{eq:shifted-afe} into the corresponding numerator moment.
For square-free $m,n$, a prime has odd exponent in $mn$ precisely
when it divides exactly one of $m,n$.  These are the prime divisors
of $k(m,n)=mn/(m,n)^2$, which is square-free.  Multiplying by
$r$ makes every exponent even if and only if the square-free part of
$r$ is $k(m,n)$, or equivalently
\[
 r=k(m,n)b^2\qquad(b\geq1).
\]
At $b=1$, the primes of $rmn$ are exactly those of $[m,n]$, so the
orthogonality main term has the factor $h([m,n])$ occurring in
Lemma~\ref{lem:mesoscopic-resonator}.
Let $\varepsilon_N=(L_2L_3/L)^{1/4}$.  The cutoff argument in
Lemma~\ref{lem:mesoscopic-resonator} shows that
\eqref{eq:mesoscopic-gal} remains valid with $\varepsilon_N$ in place
of a fixed $\varepsilon$.  Keep $b=1$ and
$k(m,n)\leq N^{\varepsilon_N}=X^{o(1)}$.  Equations
\eqref{eq:first-boundary} and \eqref{eq:mesoscopic-gal} give the
following contribution from the first half
\begin{equation}\label{eq:mesoscopic-retained}
 \left(\frac2{\pi^2}+o(1)\right)
 XH_N\mu_N^\ell\exp(G_N).
\end{equation}
For $\ell>0$, the diagonal $k=1$ has zero logarithmic weight and
at most $O_\ell(XH_N)$ actual kernel contribution. This is negligible.
For every other selected pair, $\log k\to\infty$ uniformly because
$k$ has a prime factor exceeding $\mathrm eY_0$.  Equation
\eqref{eq:first-boundary} applies with a relative $o(1)$ error.
For $\ell=0$ it applies directly also to $k=1$.
Equations \eqref{eq:dual-boundary} and
\eqref{eq:opposite-general} give the dual contribution
\begin{equation}\label{eq:mesoscopic-dual-retained}
 \left(\frac2{\pi^2}+o(1)\right)
 XH_NX^{-\alpha}(\log X)^\ell\exp(G_N^-).
\end{equation}

All other dual square terms are nonnegative by
\eqref{eq:dual-positive}. All first-half terms with $r\geq R_\ell$ are
nonnegative by \eqref{eq:first-positive}.  The remaining terms have
$r<R_\ell$.  Since every prime in $\mathcal P$ tends to infinity, this
forces $k(m,n)=1$, and hence $m=n$.  Their total possible negative
contribution is only $O_\ell(XH_N)$ by
\eqref{eq:first-finite-negative}, which is negligible compared with
\eqref{eq:mesoscopic-retained}.

We now bound the nonsquare part, including its dependence on $r$.
First restrict to $r\leq X$.  If $rmn=n_0n_1^2$, then
$\log n_0\ll LL_2/L_3$, so $g_1(n_0)=X^{o(1)}$ as in the
denominator.  The prime factors of $n_1$ outside $\mathcal P$ come
from $r$, while at most $2aL/L_3$ primes from $\mathcal P$ divide
$mn$.  We obtain
\[
 g_2(n_1)\leq\tau(r)
 \prod_{p\mid mn}(1+p^{-1/2-\delta})=X^{o(1)}
 \qquad(r\leq X),
\]
by the standard bound $\tau(r)=X^{o(1)}$ uniformly on this range and
the preceding estimate for the product.

Partial summation in $d$ bounds each orthogonality remainder by its
partial-sum bound times the supremum of the kernel and its total
variation on $[X,2X]$.  The estimates for $d\partial_d$ in
Lemma~\ref{lem:shifted-kernel} bound the total variation by the same
quantity as the kernel itself.  Since $r^{-\alpha}\leq1$ and
$X^{-\alpha}r^\alpha\leq1$ for $r\leq X$, the two kernel bounds give
\begin{align*}
 \text{nonsquare error for }r\leq X
 &\ll X^{1/2+\delta+o(1)}(\log X)^{\ell+1}
       \left(\sum_{m\in\mathcal M}f(m)\right)^2
       \sum_{r\leq X}\frac{(1+r/\sqrt X)^{-2}}{\sqrt r}\\
 &\ll X^{3/4+\delta+o(1)}(\log X)^{\ell+1}
       N\sum_{m\in\mathcal M}f(m)^2.
\end{align*}
For the last sum, the part $r\leq\sqrt X$ is $O(X^{1/4})$.
On each interval $2^v\sqrt X<r\leq2^{v+1}\sqrt X$ the bound is
$O(X^{1/4}2^{-3v/2})$.  This proves the displayed $X^{1/4}$ factor.

For $r>X$, use the trivial bound for the character sum rather than
orthogonality.  Since $\alpha=o(1)$, the global kernel estimates with
$C=10$ give
\[
 \sum_{r>X}\frac{|W^+_{\ell,\alpha,d}(r)|+
                         |W^-_{\ell,\alpha,d}(r)|}{\sqrt r}
 \ll\sum_{r>X}r^{-1/2+\alpha}\log^{\ell+1}(2Xr)
                       (r/\sqrt X)^{-10}
 \ll X^{-4}.
\]
After summing over $d,m,n$, this tail is
$O(X^{-3}N\sum_mf(m)^2)=o(XH_N)$.  Combining the two ranges yields
\begin{equation}\label{eq:mesoscopic-orth-error}
 \ll_{A_0,\ell,\delta}X^{3/4+\delta+o(1)}N
 (\log X)^{\ell+1}\sum_{m\in\mathcal M}f(m)^2=o(XH_N).
\end{equation}
Its exponent relative to $XH_N$ is
$\delta-\beta+o(1)<0$.  The numerator is bounded below by the
sum of \eqref{eq:mesoscopic-retained} and
\eqref{eq:mesoscopic-dual-retained}, up to a negligible error.
Because $R_d^2\geq0$ and the denominator is positive, the quotient of the numerator
by \eqref{eq:mesoscopic-den} is at most
$\max_{X<d\leq2X,\ d\in\Fplus}(-1)^\ell
L^{(\ell)}(1/2+\alpha,\chi_d)$.
This gives a lower bound for the required maximum.

Dividing the two main terms by
\eqref{eq:mesoscopic-den} gives
\begin{align*}
 \max_{\substack{X<d\leq2X\\d\in\Fplus}}
 (-1)^\ell L^{(\ell)}(\sigma_X,\chi_d)
 &\geq(1-o(1))\left\{
 \mu_N^\ell\exp(G_N)\right.\\
 &\hspace{26mm}\left.
 +X^{-\alpha}(\log X)^\ell\exp(G_N^-)\right\}.
\end{align*}
The orthogonality error is $o(XH_N)$, and the finite negative-square
contribution is $O_\ell(XH_N)$.
Since $\mu_N^\ell\exp(G_N)\to\infty$ uniformly for $0\leq A\leq A_0$,
both are negligible relative to the first positive main term.
They are therefore absorbed into the displayed relative $o(1)$
even when the dual contribution is small.  This proves
\eqref{eq:mesoscopic-main} with the stated uniformity.
\end{proof}

\begin{proof}[Proof of Corollary~\ref{cor:mesoscopic-exponent}]
For fixed admissible $\beta,\gamma,a$, the scales in the construction
satisfy
\begin{align*}
 \sqrt{\frac{LL_3}{L_2}}
  &=(\sqrt{1/4-\beta}+o(1))\mathcal V_X,\\
 \sqrt{LL_2L_3}
  &=(\sqrt{1/4-\beta}+o(1))\mathcal M_X.
\end{align*}
Lemma~\ref{lem:mesoscopic-resonator} gives
$G_N=(2\gamma\mathrm e^{-A}+o(1))
\sqrt{LL_3/L_2}$ and
$\mu_N=(2\gamma\mathrm e^{-A}+o(1))\sqrt{LL_2L_3}$.
The factor $\mu_N^\ell$ has logarithm $O_\ell(\log_2X)$
uniformly for bounded $A$, so using just the first term of
\eqref{eq:mesoscopic-main} gives the leading coefficient
$2\gamma\sqrt{1/4-\beta}\,\mathrm e^{-A}$.

To let this coefficient approach $\mathrm e^{-A}$ without assuming
uniformity in $\beta,\gamma,a$, take, for integers $j\geq1$,
\[
 \beta_j=\frac1{8j},\qquad
 \gamma_j=1-\frac1{j+1},\qquad
 a_j=1+\frac1{2(j+1)}.
\]
These satisfy $a_j\gamma_j<1$ and converge to $0,1,1$, respectively.
For each fixed $j$, all the estimates above hold uniformly for
$0\leq A\leq A_0$.  Choose increasing thresholds $X_j\to\infty$
so that the relative errors and the errors in the normalized
free-energy and mean estimates are at most $1/j$ for $X\geq X_j$.
Use $\delta=\beta_j/2$ in the arithmetic estimates.
For $X_j\leq X<X_{j+1}$ use the parameters indexed by $j$.
Then $j\to\infty$ with $X$, all these errors tend to zero, and
$2\gamma_j\sqrt{1/4-\beta_j}\to1$.
This proves the corollary uniformly for $0\leq A\leq A_0$.

For completeness, the full two-term expression at this leading
precision has the same exponential coefficient.  After the diagonal
choice, the logarithm of its dual term is
\[
 -\frac{A\log X}{\log_2X}
 +(\mathrm e^A+o(1))\mathcal V_X+\ell\log_2X.
\]
On $[0,A_0]$,
$\mathrm e^A-\mathrm e^{-A}\leq2\cosh(A_0)A$, whereas
\[
 \frac{\log X/\log_2X}{\mathcal V_X}
 =\sqrt{\frac{\log X}{\log_2X\log_3X}}\longrightarrow\infty.
\]
Consequently
$-A\log X/\log_2X+
(\mathrm e^A-\mathrm e^{-A})\mathcal V_X\leq0$
for every $0\leq A\leq A_0$ once $X$ is large.
The logarithms of both derivative weights are
$O_\ell(\log_2X)=o(\mathcal V_X)$, and taking the sum adds at most
$\log2$ to the larger logarithm.  The two-term lower bound thus has
logarithm $(\mathrm e^{-A}+o(1))\mathcal V_X$ uniformly.
This comparison keeps the conductor penalty
$-A\log X/\log_2X$ explicit.  For fixed $A>0$ it is larger in magnitude
than $\mathcal V_X$ and makes the dual term negligible.
\end{proof}

\begin{proof}[Proof of Theorem~\ref{thm:boundary-layer}]
Fix $0<\beta<1/4$, $0<\gamma<1$, and $a>1$ with $a\gamma<1$, and put
\[
 N=\lfloor X^{1/4-\beta}\rfloor,\qquad
 \alpha=\frac A{\log_2X}.
\]
Use the resonator supplied by Lemma~\ref{lem:mesoscopic-resonator}.
The hypothesis $A\mathcal V_X\to0$ implies
\[
 \alpha\sqrt{LL_2L_3}
 =(\sqrt{1/4-\beta}+o(1))A\mathcal V_X\to0,
\]
which is the relative-error hypothesis of
Corollary~\ref{cor:opposite-tilts}.  We also have
$N^{\varepsilon_N}=X^{o(1)}$.

The denominator is again \eqref{eq:mesoscopic-den}.  In the numerator,
write every square index as $r=k(m,n)b^2$.  Keep $b=1$ and
$k(m,n)\leq N^{\varepsilon_N}$.  Apart from $k=1$, which has only
baseline mass, every nontrivial $k(m,n)$ is larger than $R_\ell$ for
large $X$.  Equations \eqref{eq:first-boundary} and
\eqref{eq:boundary-positive-tilt} give the first contribution
\begin{equation}\label{eq:boundary-first-contribution}
 \left(\frac2{\pi^2}+o(1)\right)XH_N
 \mu_N^\ell\exp(G_N^0).
\end{equation}
Equations \eqref{eq:dual-boundary} and
\eqref{eq:boundary-negative-tilt} give the dual contribution
\begin{equation}\label{eq:boundary-dual-contribution}
 \left(\frac2{\pi^2}+o(1)\right)XH_N
 X^{-\alpha}(\log X)^\ell\exp(G_N^0).
\end{equation}

It remains to justify that the other square terms cannot cancel
these two positive quantities.  By \eqref{eq:dual-positive}, every
dual term is nonnegative.  By \eqref{eq:first-positive}, the same is
true of the first half when $r\geq R_\ell$.  If $r<R_\ell$, then
$r=k(m,n)b^2$ forces $k(m,n)=1$ for large $X$, because every prime in
$\mathcal P$ exceeds $R_\ell$.  Thus $m=n$, and
\eqref{eq:first-finite-negative} shows that the total possible negative
square contribution is only $O_\ell(XH_N)$.  This is negligible compared
with \eqref{eq:boundary-first-contribution}, since $G_N^0\to\infty$.

The nonsquare orthogonality error is bounded as in
\eqref{eq:mesoscopic-orth-error}, now without using any negative power
of $X^\alpha$:
\[
 \ll_{\ell,\delta}X^{3/4+\delta+o(1)}N
 (\log X)^{\ell+1}\sum_{m\in\mathcal M}f(m)^2=o(XH_N)
\]
after $0<\delta<\beta$ is fixed.  Dividing the sum of
\eqref{eq:boundary-first-contribution} and
\eqref{eq:boundary-dual-contribution} by
\eqref{eq:mesoscopic-den} proves \eqref{eq:boundary-precise}.
Because $R_d^2\geq0$, this quotient is a weighted average bounded
above by the maximum in \eqref{eq:boundary-precise}.  The first
term alone tends to infinity after division by $XH_N$, so the additive
errors also give a relative $o(1)$ for the sum of both terms.
All the estimates are uniform on a prescribed range
$0\leq A\leq a_X$ with $a_X\mathcal V_X\to0$.

For the additional diagonal assertion, use the parameter sequences
from the proof of Corollary~\ref{cor:mesoscopic-exponent}.
For each fixed $j$, increase the threshold $X_j$ so that the relative
errors in Corollary~\ref{cor:opposite-tilts} and
\eqref{eq:boundary-precise}, and the errors in the normalized mean
and free-energy estimates, are at most $1/j$ for $X\geq X_j$,
uniformly for $0\leq A\leq a_X$.
Such thresholds exist because the actual-shift condition of that
corollary holds uniformly on the prescribed range.
For $X_j\leq X<X_{j+1}$, use $\beta_j,\gamma_j,a_j$ and
$\delta=\beta_j/2$.  The fixed-parameter estimates
\[
 \mu_N=(2\gamma_j+o(1))\sqrt{LL_2L_3},\qquad
 G_N^0=(2\gamma_j+o(1))\sqrt{\frac{LL_3}{L_2}}
\]
then combine with $2\gamma_j\sqrt{1/4-\beta_j}\to1$ to give
$\mu_N/\mathcal M_X\to1$ and
$G_N^0=(1+o(1))\mathcal V_X$, uniformly on that range.
The finite products in \eqref{eq:boundary-precise} are used without
approximation throughout this choice, so its relative estimate also
persists.

For these choices, the ratio of the dual contribution to the
first contribution is
\[
 (1+o(1))\frac{X^{-\alpha}(\log X)^\ell}{\mathcal M_X^\ell}.
\]
Since
$\log\mathcal M_X=
\tfrac12(\log_2X+\log_3X+\log_4X)$,
the logarithm of the ratio is
\[
 -\frac{A\log X}{\log_2X}
 +\frac\ell2(\log_2X-\log_3X-\log_4X)+o(1).
\]
For $\ell\geq1$ this equals
$-(A-A_{\ell,*})\log X/\log_2X+o(1)$, proving all three
comparisons following \eqref{eq:boundary-location}.
Also $A_{\ell,*}=O_\ell((\log_2X)^2/\log X)$, and hence
\[
 A_{\ell,*}\mathcal V_X
 \ll_\ell\frac{(\log_2X)^{3/2}\sqrt{\log_3X}}{\sqrt{\log X}}
 \longrightarrow0.
\]
For $\ell=0$ the ratio is $(1+o(1))X^{-\alpha}$ and the formal
equality is at $A=0$. These calculations compare the two contributions.
See Remark~\ref{rem:sharpness} for the scope of this comparison.
\end{proof}

\begin{proof}[Proof of Corollary~\ref{cor:boundary-exponent}]
Use the admissible diagonal choices in
Theorem~\ref{thm:boundary-layer}.  They give
$G_N^0=(1+o(1))\mathcal V_X$ and
$\mu_N=(1+o(1))\mathcal M_X$, while
$0<X^{-\alpha}\leq1$.
For $\ell\geq1$, the factor in braces in
\eqref{eq:boundary-precise} is therefore bounded below by
$(1+o(1))\mathcal M_X^\ell$ and above by
$(1+o(1))\mathcal M_X^\ell+(\log X)^\ell$.
Its logarithm is $O_\ell(\log_2X)=o(\mathcal V_X)$.
For $\ell=0$, the factor lies between $1$ and $2$ and its logarithm
is bounded.  Taking logarithms in \eqref{eq:boundary-precise}
proves \eqref{eq:boundary-layer}, with the stated uniformity.
\end{proof}

\begin{remark}
\label{rem:mesoscopic-optimality}
One can optimize the common scalar amplitude while keeping the block
profile and the type of occupancy cutoffs fixed.  If the function in
\eqref{eq:mesoscopic-f} is multiplied by $\lambda>0$, the typical block
occupancy and hence the entropy cost are multiplied by $\lambda^2$.
The cardinality condition becomes
$\lambda^2\gamma\leq1+o(1)$, whereas the free-energy coefficient before
the substitution $N=X^{1/4-o(1)}$ is
$2\lambda\gamma\mathrm e^{-A}$.  Thus
\[
 \lambda\gamma\leq\sqrt\gamma\leq1.
\]
The limiting choice $\lambda=1$, $\gamma\uparrow1$ used above is
therefore optimal within this common-rescaling family.  This computation
does not optimize over different block profiles, endpoints, or cutoff
rules.  It makes no general optimality claim for the resonance method.
\end{remark}

\section{Preparations for the proof of Theorem \ref{thm2}}\label{sec5}

Fix $0<b<1$, choose an admissible $r$ in \eqref{eq:Delta-r}, and take
\[
 \kappa>\kappa_r(\sigma),\qquad
 0<a<\frac{1-2\kappa(1-\sigma)}{\alpha(b)},
 \qquad Y=a\log X\log_2X,
\]
where
\[
 \alpha(b)=2\log\frac{(1+b)^2}{1+b^2}.
\]
These choices ensure that $1-2\kappa(1-\sigma)>0$.
Let $v$ be completely multiplicative with $v(p)=b$ for $p\leq Y$ and
$v(p)=0$ otherwise, and set
\[
 R_d=\sum_{n\geq1}v(n)\chi_d(n)
 =\prod_{p\leq Y}(1-b\chi_d(p))^{-1}.
\]
Every local factor is positive because $\chi_d(p)\in\{-1,0,1\}$ and
$0<b<1$.  We can therefore use $R_d^2$ as a nonnegative weight
in the resonance average.

Define
\begin{align}\label{eq:Dproduct}
 \mathcal D&=\sum_{mn=\square}v(m)v(n)h(mn)\\
 &=\prod_{p\leq Y}\left(1-h(p)+h(p)
 \frac{1+b^2}{(1-b^2)^2}\right)
 =X^{a\vartheta(b)+o(1)},\nonumber
\end{align}
where
$\vartheta(b)=\log((1+b^2)/(1-b^2)^2)$.
This is the constant-coefficient Euler-factor calculation in
\cite[arXiv version, Section~6, (30)]{DarbarMaiti}, with the local
density $h(p)$ used here.  The density at $p=2$ changes
only a fixed local factor, so the exponent of $X$ is unchanged.

For each $p\leq Y$, define
\[
 D_p=1-h(p)+h(p)\frac{1+b^2}{(1-b^2)^2},\qquad
 c_p=\frac{2b\,h(p)}{(1-b^2)^2D_p}.
\]
The finite product and its logarithmic derivative are
\[
 Z_Y(\sigma)=\prod_{p\leq Y}(1+c_pp^{-\sigma}),\qquad
 \mu_Y=\sum_{p\leq Y}
 \frac{c_pp^{-\sigma}}{1+c_pp^{-\sigma}}\log p.
\]
These quantities correspond to the fixed choices of
$\sigma,b,a,r,\kappa$.

\begin{proposition}\label{prop:fixed-transfer}
With the choices above,
\begin{align}\label{eq:fixed-den}
 \sum_{\substack{X<d\leq2X\\d\in\Fplus}}R_d^2
 &=\left(\frac2{\pi^2}+o(1)\right)X\mathcal D,
\end{align}
and
\begin{equation}\label{eq:fixed-precise}
 \frac{\displaystyle\sum_{\substack{X<d\leq2X\\d\in\Fplus}}
 (-1)^\ell L^{(\ell)}(\sigma,\chi_d)R_d^2}
 {\displaystyle\sum_{\substack{X<d\leq2X\\d\in\Fplus}}R_d^2}
 \geq(1-o(1))\mu_Y^\ell Z_Y(\sigma).
\end{equation}
At the leading exponential scale, this gives
\begin{align}\label{eq:fixed-num}
 &\frac{\displaystyle
 \sum_{\substack{X<d\leq2X\\d\in\Fplus}}
 (-1)^\ell L^{(\ell)}(\sigma,\chi_d)R_d^2}
 {\displaystyle
 \sum_{\substack{X<d\leq2X\\d\in\Fplus}}R_d^2}\\
 &\quad\geq
 \exp\!\left(\left(\frac{2b}{1+b^2}
 \frac{a^{1-\sigma}}{1-\sigma}+o(1)\right)
 \frac{(\log X)^{1-\sigma}}{(\log_2X)^\sigma}\right).
 \nonumber
\end{align}
\end{proposition}

\begin{proof}
Apply the truncation estimate underlying Lemma~\ref{lem:burgess-trunc}
with the common cutoff $K=X^\kappa$.  Uniformly for $X<d\leq2X$,
\[
 (-1)^\ell L^{(\ell)}(\sigma,\chi_d)
 =\sum_{k\leq K}\frac{\chi_d(k)(\log k)^\ell}{k^\sigma}+o(1).
\]
Here $d\asymp X$, so the strict inequality
$\kappa>\kappa_r(\sigma)$ gives the same power saving for this cutoff.
The uniform error contributes $o(\sum_dR_d^2)$ to the numerator.

We use the multiplicative-majorant calculation in
\cite[arXiv version, Section~6, the estimates preceding (30)]{DarbarMaiti},
with the orthogonality estimate in Lemma~\ref{lem:orth} for
the present subfamily.  This gives
\[
 \sum_{m\geq1}v(m)\prod_{p\mid m}(1+p^{-1/2-\eps})
 =X^{-a\log(1-b)+o(1)},\qquad
 g_1((kmn)_0)=X^{o(1)}
\]
for $k\leq K$ and $v(m)v(n)\ne0$.
For the second estimate, the primes from $mn$ lie below $Y$,
while those from $k$ have product at most $K$.  The bound for
$g_1$ in Lemma~\ref{lem:orth} is uniformly $X^{o(1)}$.
The exponents in $m,n$ do not enter its argument.
The factor from $k$ is $O_\eps(k^\eps)$, and
\[
 \sum_{k\leq K}\frac{(\log k)^\ell}{k^{\sigma-\eps}}
 \leq X^{\kappa(1-\sigma)+O(\eps)+o(1)}.
\]
The denominator error is
\begin{equation}\label{eq:fixed-den-error}
 \ll X^{1/2-2a\log(1-b)+O(\eps)+o(1)},
\end{equation}
while the numerator error, including the additional $k$-sum, is
\begin{equation}\label{eq:fixed-num-error}
 \ll X^{1/2+\kappa(1-\sigma)-2a\log(1-b)+O(\eps)+o(1)}.
\end{equation}
First expand the Euler products with prime exponents at most
$J$, then let $J\to\infty$.  The convergent majorant above
bounds the errors absolutely.  The square main terms have
nonnegative coefficients.

The square denominator main term is $2X\mathcal D/\pi^2$, with
$\mathcal D=X^{a\vartheta(b)+o(1)}$ as in \eqref{eq:Dproduct}.
The corresponding local-factor calculation also appears in
\cite[arXiv version, (30)]{DarbarMaiti}.
Since
\[
 -2\log(1-b)-\vartheta(b)=\frac{\alpha(b)}2,\qquad
 -\frac12+\kappa(1-\sigma)+\frac{a\alpha(b)}2<0,
\]
choosing $\eps>0$ sufficiently small after the fixed parameters
makes both \eqref{eq:fixed-den-error} and
\eqref{eq:fixed-num-error} equal to $o(X\mathcal D)$.
This proves \eqref{eq:fixed-den}.

In the numerator, restrict $k$ to square-free integers supported on $p\leq Y$.
At a prime dividing $k$, the condition $kmn=\square$ requires the
two resonator exponents to have odd sum.  Their local
contribution is
\[
 h(p)\sum_{\substack{e,f\geq0\\e+f\ \mathrm{odd}}}b^{e+f}
 =\frac{2b\,h(p)}{(1-b^2)^2}.
\]
Division by the denominator factor $D_p$ gives $c_p$.
At primes not dividing $k$, the factor is $D_p$.
The resulting square contribution, divided by
$2X\mathcal D/\pi^2$, is
\begin{equation}\label{eq:fixed-square-quotient}
 \sum_{\substack{k\leq K,\ k\ \mathrm{square\mbox{-}free}\\
                   P^+(k)\leq Y}}
 \frac{c_k(\log k)^\ell}{k^\sigma},
 \qquad c_k=\prod_{p\mid k}c_p.
\end{equation}
All omitted square terms are nonnegative.
Moreover, $c_p=2b/(1+b^2)+O_b(p^{-1})$, so
Lemma~\ref{lem:free-energy} evaluates the unrestricted sum with
relative precision as $(1+o(1))\mu_Y^\ell Z_Y(\sigma)$.

We next check that the cutoff preserves the derivative weight.
Under the Bernoulli measure of Lemma~\ref{lem:free-energy},
division of \eqref{eq:fixed-square-quotient} by $Z_Y(\sigma)$ gives
$\mathbb E(S_Y^\ell\ind_{S_Y\leq\log K})$.
Exponential Markov with exponent $1/\log Y$ yields
\begin{align*}
 \mathbb P(S_Y>\log K)
 &\leq K^{-1/\log Y}\prod_{p\leq Y}
 \left(1+\frac{c_pp^{-\sigma}}{1+c_pp^{-\sigma}}
                  (p^{1/\log Y}-1)\right)\\
 &\leq\exp\left(-\frac{\kappa\log X}{\log Y}
 +O_{\sigma,b}\left(\frac{Y^{1-\sigma}}{\log Y}\right)\right)
 \leq\exp\left(-\frac{\kappa\log X}{2\log Y}\right).
\end{align*}
Here $p^{1/\log Y}\leq e$ and $Y^{1-\sigma}=o(\log X)$.
Since $S_Y\leq\sum_{p\leq Y}\log p\ll Y$,
\[
 \mathbb E(S_Y^\ell\ind_{S_Y>\log K})
 \ll_\ell Y^\ell
 \exp\left(-\frac{\kappa\log X}{2\log Y}\right)
 =o(\mu_Y^\ell).
\]
For $\ell=0$ this is the probability estimate itself.
Equation~\eqref{eq:fixed-square-quotient} equals
$(1+o(1))\mu_Y^\ell Z_Y(\sigma)$.
Together with the denominator and error estimates, this proves
\eqref{eq:fixed-precise}.

Equation~\eqref{eq:mu-asymp} evaluates $\mu_Y$.
For the coarse consequence, first-order expansion of the local
logarithms and the prime number theorem give
\[
 \log Z_Y(\sigma)
 =\sum_{p\leq Y}c_pp^{-\sigma}+O_{\sigma,b}(1)
 =\left(\frac{2b}{1+b^2}\frac1{1-\sigma}+o(1)\right)
 \frac{Y^{1-\sigma}}{\log Y}.
\]
The $O(1)$ error follows from convergence of
$\sum_p p^{-2\sigma}$ and $\sum_p p^{-1-\sigma}$.
Substitute $Y=a\log X\log_2X$ and $\log Y\sim\log_2X$, and absorb
$\ell\log\mu_Y=O_{\ell,\sigma,b,a}(\log_2X)$ into the coarse exponent
error.  This proves \eqref{eq:fixed-num}.
\end{proof}

\section{Proof of Theorem \ref{thm2}}\label{sec6}

\begin{proof}
Fix any admissible parameter choice in Theorem~\ref{thm2}.
The coefficient $c_p$ in its statement agrees with the
coefficient in Proposition~\ref{prop:fixed-transfer}.
Multiplying the local factor $D_p$ there by $(1-b^2)^2$ gives the
denominator in the theorem's definition of $c_p$.
The quantities $Z_Y(\sigma)$ and $\mu_Y$ are therefore the product and sum
in \eqref{eq:fixed-precise}.

The denominator in that weighted inequality is positive for all
sufficiently large $X$, by \eqref{eq:fixed-den}.  Since $R_d^2\geq0$,
\[
 \frac{\displaystyle\sum_{\substack{X<d\leq2X\\d\in\Fplus}}
 (-1)^\ell L^{(\ell)}(\sigma,\chi_d)R_d^2}
 {\displaystyle\sum_{\substack{X<d\leq2X\\d\in\Fplus}}R_d^2}
 \leq
 \max_{\substack{X<d\leq2X\\d\in\Fplus}}
 (-1)^\ell L^{(\ell)}(\sigma,\chi_d).
\]
Combining this comparison with \eqref{eq:fixed-precise} proves
\eqref{eq:fixed-main-relative}.  The comparison preserves both the
relative error and the finite product.
\end{proof}

\begin{proof}[Proof of Corollary~\ref{cor:fixed-exponent}]
We evaluate and optimize the exponential size of the bound just proved.
For every fixed admissible choice, the evaluations in Proposition
\ref{prop:fixed-transfer} give
\[
 \log Z_Y(\sigma)=
 \left(\frac{2b}{1+b^2}\frac{a^{1-\sigma}}{1-\sigma}+o(1)\right)
 \frac{(\log X)^{1-\sigma}}{(\log_2X)^\sigma}.
\]
They also give
\[
 \log\mu_Y=(1-\sigma)(\log_2X+\log_3X)+O(1).
\]
Since
\[
 \frac{\log_2X}
 {(\log X)^{1-\sigma}/(\log_2X)^\sigma}
 =\frac{(\log_2X)^{1+\sigma}}{(\log X)^{1-\sigma}}
 \longrightarrow0,
\]
the logarithms of $\mu_Y^\ell$ and the relative factor
$1-o(1)$ are absorbed by the unquantified error in the
exponential evaluation.  Theorem~\ref{thm2} gives
\[
 \max_{\substack{X<d\leq2X\\d\in\Fplus}}
 (-1)^\ell L^{(\ell)}(\sigma,\chi_d)
 \geq
 \exp\!\left(
 \left(\frac{2b}{1+b^2}\frac{a^{1-\sigma}}{1-\sigma}+o(1)\right)
 \frac{(\log X)^{1-\sigma}}{(\log_2X)^\sigma}\right).
\]

Choose $r$ attaining $\Delta(\sigma)$.
For fixed $b$, take $\kappa>\kappa_r(\sigma)$ arbitrarily close to
$\kappa_r(\sigma)$.  Then take $a$ arbitrarily close from below
to $(1-2\kappa(1-\sigma))/\alpha(b)$.
The limiting value of $a$ is
\[
 \frac{\Delta_r(\sigma)}{\alpha(b)}.
\]
Finally let $b\uparrow1$.
The limits $2b/(1+b^2)\to1$ and $\alpha(b)\to2\log2$ show that the
coefficient in the last exponential approaches
\[
 \frac1{1-\sigma}
 \left(\frac{\Delta(\sigma)}{2\log2}\right)^{1-\sigma}
 =C_{\rm B}(\sigma).
\]

Each limit is taken through fixed parameter choices.
No estimate uniform as $b\uparrow1$ is needed.
To obtain one asymptotic statement, take a sequence of admissible fixed
choices whose coefficients tend to $C_{\rm B}(\sigma)$.
For its $j$-th member choose a threshold beyond which the error in the
normalized logarithmic bound is at most $1/j$.
Make the thresholds strictly increasing and tending to infinity, and
use the $j$-th fixed choice between its threshold and the next one.
Then the parameter loss and the normalized logarithmic error both
tend to zero as $X\to\infty$.
This gives the coefficient $C_{\rm B}(\sigma)+o(1)$ in
Corollary~\ref{cor:fixed-exponent}.
\end{proof}

\section{Preparations for the proof of Theorem \ref{thm3}}\label{sec7}

Put
\[
 y=\frac{\log X}{100\log_2X},\qquad
 Q=\prod_{p\leq y}p,
\]
and define the saturated divisor resonator
\[
 R_d=\sum_{m\mid Q}\chi_d(Qm).
\]
The extra factor $Q$ ensures that all selected square terms have the
same local density $h(Q)$.

\begin{proposition}\label{prop:saturated}
Let $\ell\geq0$ be fixed, $1/2<\sigma\leq1$, and suppose
\[
 0\leq(1-\sigma)\log_2X\leq\frac12\log_3X.
\]
Then, uniformly in this range,
\begin{equation}\label{eq:saturated}
 \max_{\substack{X<d\leq2X\\d\in\Fplus}}
 (-1)^\ell L^{(\ell)}(\sigma,\chi_d)
 \geq(1-o(1))
\sum_{k\in\Ssmooth(y)}\frac{(\log k)^\ell}{k^\sigma}.
\end{equation}
\end{proposition}

\begin{proof}
Fix a small $\eta>0$ and take $K=X^{1/(2\sigma)+\eta}$.
The P\'olya--Vinogradov partial-summation bound recorded after
Lemma~\ref{lem:burgess-trunc}, with the standard input
\cite[Chapter~12]{IK}, gives
\[
 \left|\sum_{k>K}\frac{\chi_d(k)(\log k)^\ell}{k^\sigma}\right|
 \ll_\ell\sqrt X\log X\,K^{-\sigma}(\log K)^\ell
 \ll_\ell X^{-\eta\sigma}(\log X)^{\ell+1}=o(1).
\]
The bound is uniform for $1/2<\sigma\leq1$, since the
partial-summation integral has a bounded factor $1/\sigma$.

There are $2^{\pi(y)}$ divisors of $Q$, and the prime number theorem
gives
\[
 2^{\pi(y)}=X^{o(1)},\qquad Q^2=X^{o(1)}.
\]
In the denominator expansion, $Q^2mn$ is a square, with $m,n\mid Q$,
exactly when $m=n$.  The $2^{\pi(y)}$ square
pairs each have density $h(Q)$.  Lemma~\ref{lem:orth} gives
\[
 \sum_dR_d^2=\left(\frac2{\pi^2}+o(1)\right)
 Xh(Q)2^{\pi(y)}.
\]
The sum here and below is over $X<d\leq2X$, $d\in\Fplus$.
The error is negligible because all the square-free parts are
$X^{o(1)}$, there are $2^{2\pi(y)}=X^{o(1)}$ pairs, and
$h(Q)\asymp1/\log y=X^{-o(1)}$.

In the numerator, consider only $y$-smooth integers $k$.
At each prime $p\leq y$, the exponents of $m$ and $n$ are in
$\{0,1\}$.
If the exponent of $p$ in $k$ is even, their allowed pairs are
$(0,0)$ and $(1,1)$.  If it is odd, the allowed pairs are
$(0,1)$ and $(1,0)$.
The condition $kQ^2mn=\square$ has exactly $2^{\pi(y)}$ ordered solutions,
independently of $k$.  The factor $Q^2$ ensures that the prime support
of every such square contains all primes at most $y$, so
$h(kQ^2mn)=h(Q)$.
The corresponding numerator main term is
\[
 \frac2{\pi^2}Xh(Q)2^{\pi(y)}
 \sum_{\substack{k\leq K\\P^+(k)\leq y}}
 \frac{(\log k)^\ell}{k^\sigma}.
\]
All other square terms are nonnegative and can be discarded for the
lower bound.

The total character-sum error in the numerator is
\[
 \ll X^{1/2+O(\eps)+o(1)}2^{2\pi(y)}
 \sum_{k\leq K}\frac{(\log k)^\ell}{k^{\sigma-\eps}}.
\]
In this range $1-\sigma=o(1)$ uniformly and
$\log K/\log X=1/(2\sigma)+\eta=1/2+\eta+o(1)$.
Integral comparison bounds the $k$-sum by
$K^{1-\sigma+\eps}(\log K)^{\ell+1}$.
A sufficiently small fixed $\eps>0$ makes the
whole error $O(X^{3/4})$, uniformly in $\sigma$.
This is $o(Xh(Q)2^{\pi(y)})$.
The uniform $o(1)$ truncation error from the first paragraph contributes
$o(\sum_dR_d^2)$ as well.

Lemma~\ref{lem:friable} removes the cutoff $k\leq K$ with relative
error $o(1)$.  To make the uniformity explicit, use a cutoff independent of $\sigma$.
First restrict to $k\leq X^{1/2}$, which lies below $K$.
The lemma's truncation assertion gives the full sum with the
required relative error for this smaller cutoff.
Dividing the numerator by the denominator yields the full smooth sum
with factor $1-o(1)$.  Finally, the maximum is at least this average
because the weights $R_d^2$ are nonnegative.  This proves
\eqref{eq:saturated}.
\end{proof}

\section{Proof of Theorem \ref{thm3}}\label{sec8}

\begin{proof}[Proof of Theorem~\ref{thm3}]
Set $\sigma=1-A/\log_2X$ with $A\geq0$ fixed.
For the choice of $y$ in Section~\ref{sec7},
\[
 \log y=\log_2X-\log_3X-\log100,
 \qquad (1-\sigma)\log y=A+o(1).
\]
Lemma~\ref{lem:friable}, uniformly on a compact interval containing
these parameters, gives
\[
 \sum_{k\in\Ssmooth(y)}\frac{(\log k)^\ell}{k^\sigma}
 =\bigl(D_\ell(A)+o(1)\bigr)(\log_2X)^{\ell+1}.
\]
Continuity of $D_\ell$ allows us to replace
$(1-\sigma)\log y$ by $A$.
The relation $\log y/\log_2X\to1$ allows us to replace the logarithmic factor.
Since the hypotheses of Proposition~\ref{prop:saturated} hold for every
fixed $A$ once $X$ is large, insertion of this evaluation into
\eqref{eq:saturated} proves Theorem~\ref{thm3}.
\end{proof}

\begin{proof}[Proof of Corollary~\ref{thm4}]
Taking $A=0$ in Theorem~\ref{thm3} gives $\sigma_A=1$ and the
constant $D_\ell(0)$, which is the stated conclusion.
\end{proof}

\section{Proof of Theorem \ref{thm5}}\label{sec10}

For $1\leq A\leq\tfrac12\log_3X$, Proposition~\ref{prop:saturated}
already provides the required transfer uniformly.
Farther from the $1$-line, we use the two-piece approximate functional
equation.  A fixed-$\eps$ orthogonality estimate then suffices
throughout the range.

\begin{proposition}\label{prop:shifted-saturated}
Let $1/2<\sigma<1$, put
\[
 A=(1-\sigma)\log_2X,\qquad
 B=(\sigma-\tfrac12)\log_2X,
\]
and suppose that $A\geq\tfrac12\log_3X$ and $B\to\infty$.
Uniformly in this range,
\begin{equation}\label{eq:shifted-saturated}
 \max_{\substack{X<d\leq2X\\d\in\Fplus}}
 (-1)^\ell L^{(\ell)}(\sigma,\chi_d)
 \geq(1-o(1))
 \sum_{n\in\Ssmooth(\log X)}\frac{(\log n)^\ell}{n^\sigma}.
\end{equation}
\end{proposition}

\begin{proof}
Set $y=\log X$, $Q=\prod_{p\leq y}p$, and
\[
 R_d=\sum_{m\mid Q}\chi_d(Qm).
\]
There are $2^{\pi(y)}=X^{o(1)}$ divisors in this sum.
Although $Q=X^{1+o(1)}$, the square-free part of $Q^2mn$ has
logarithm $O(\log X)$.  Its $g_1$ factor is therefore $X^{o(1)}$.
The other factor in Lemma~\ref{lem:orth} is bounded by
\[
 \prod_{p\leq y}(1+p^{-1/2-\eps})=X^{o(1)}.
\]
As in the preceding saturated calculation, only $m=n$ gives a square
in the denominator.  Each of the $2^{\pi(y)}$ square pairs has density
$h(Q)$, giving
\begin{equation}\label{eq:shifted-saturated-den}
 \sum_dR_d^2=\left(\frac2{\pi^2}+o(1)\right)
 Xh(Q)2^{\pi(y)}.
\end{equation}
The relative error is $o(1)$, since the absolute error is
$X^{1/2+\eps+o(1)}2^{2\pi(y)}$ and $h(Q)\asymp1/\log y$.

Insert \eqref{eq:shifted-afe} with $\alpha=\sigma-1/2$ into the
numerator moment.
For every $y$-smooth $r$, the condition $rQ^2mn=\square$, with
$m,n\mid Q$, has exactly $2^{\pi(y)}$ ordered solutions.
This follows from the prime-by-prime parity count in
Proposition~\ref{prop:saturated}.
Moreover, $h(rQ^2mn)=h(Q)$ for each solution.
Fix $0<\eta<1/8$ and select the first-half square terms with
$r\leq X^\eta$.
Equation~\eqref{eq:shifted-first-small} gives their contribution as
\begin{equation}\label{eq:shifted-saturated-main}
 \left(\frac2{\pi^2}+o(1)\right)Xh(Q)2^{\pi(y)}
 \sum_{\substack{r\leq X^\eta\\P^+(r)\leq y}}
 \frac{(\log r)^\ell}{r^\sigma}.
\end{equation}
To see that the small-kernel error is relative $o(1)$, sum its bound
$X^{-c_\eta}r^{-\sigma}(\log X)^\ell$.
For $r\geq2$, the unweighted term is at most
$(\log2)^{-\ell}$ times the weighted term.  The single term $r=1$
contributes only $O(X^{-c_\eta}(\log X)^\ell)$ after normalization.
The selected weighted sum is bounded below by a positive constant
using $r=2$, or $r=1$ when $\ell=0$.
The error in \eqref{eq:shifted-saturated-main} is relative
$o(1)$ uniformly in the stated range.

Every dual square term is nonnegative by \eqref{eq:dual-positive}.
The first-half square terms with $r\geq R_\ell$ are nonnegative by
\eqref{eq:first-positive}.  Only the finite set $r<R_\ell$ may give a
negative contribution, and \eqref{eq:first-finite-negative} bounds its
total by $O_\ell(Xh(Q)2^{\pi(y)})$.

We bound the total orthogonality remainder, including the remainders
for square indices.  Apply Lemma~\ref{lem:orth} to the kernels by
partial summation in $d$.
The bounds for $d\partial_dW^\pm$ have the same size as those for the
undifferentiated kernels.
The factors involving $m,n$ contribute $X^{o(1)}$ per pair, and those
involving $r$ can be absorbed into $r^{O(\eps)}$.
For the first kernel, the sum over $r$ is bounded by
\[
 \sum_{r\geq1}r^{-\sigma+O(\eps)}
 \log^\ell(2Xr)(1+r/\sqrt X)^{-C}
 \ll_\ell X^{(1-\sigma)/2+O(\eps)}(\log X)^{\ell+1}.
\]
For the dual kernel, the corresponding bound follows from
\[
 X^{-\alpha}\sum_{r\geq1}r^{-1/2+\alpha+O(\eps)}
 \log^{\ell+1}(2Xr)(1+r/\sqrt X)^{-C}
 \ll_\ell X^{1/4-\alpha/2+O(\eps)}(\log X)^{\ell+1}.
\]
Choose $C$ sufficiently large so that the tails converge.
Since $1/4-\alpha/2=(1-\sigma)/2\leq1/4$, multiplication by the
orthogonality factor gives, after reducing its fixed parameter if
necessary, a total error
\[
 \ll_{\ell,\eps}X^{3/4+\eps+o(1)}2^{2\pi(y)}(\log X)^{\ell+1}
 =o(Xh(Q)2^{\pi(y)})
\]
for any fixed sufficiently small $\eps<1/8$.

It remains to evaluate the selected smooth sum and remove its cutoff.
Define the product and mean
\[
 Z_y(\sigma)=\prod_{p\leq y}(1-p^{-\sigma})^{-1},\qquad
 \mu_y=\sum_{p\leq y}\frac{\log p}{p^\sigma-1}.
\]
After division by $Z_y(\sigma)$, the terms of its Dirichlet series
define a probability distribution on the smooth integers $n$.
Their prime exponents are independent, with probabilities
\[
 \mathbb P(v_p(n)=e)=(1-p^{-\sigma})p^{-e\sigma},
 \qquad e=0,1,2,\ldots.
\]
The mean of $\log n$ is $\mu_y$.
Expanding each local denominator gives
\[
 \mu_y=\sum_{p\leq y}\frac{\log p}{p^\sigma}
 +O\left(\sum_{p\leq y}\frac{\log p}{p^{2\sigma}}\right).
\]
The error is $O(\log y)$ uniformly for $\sigma>1/2$.
The prime number theorem and partial summation evaluate the first sum
uniformly as $(1+o(1))y^{1-\sigma}/(1-\sigma)$ in the present range.
The lower-endpoint contribution is negligible since $A\to\infty$.
We obtain
\begin{equation}\label{eq:shifted-saturated-mean}
 \mu_y=(1+o(1))\frac{y^{1-\sigma}}{1-\sigma}
       =(1+o(1))\frac{\e^A\log_2X}{A}.
\end{equation}
For every fixed $j\geq2$, the geometric local factors give
\[
 \kappa_j(\log n)
 =\sum_{p\leq y}(\log p)^j
       \sum_{e\geq1}e^{j-1}p^{-e\sigma}
 \ll_j\mu_y(\log y)^{j-1}.
\]
The bound is uniform since $p^{-\sigma}\leq2^{-1/2}$.
By \eqref{eq:shifted-saturated-mean},
$\log y/\mu_y=(1+o(1))A\mathrm e^{-A}\to0$ uniformly.
The normalized-cumulant argument in the proof of
Lemma~\ref{lem:free-energy} now yields
\begin{equation}\label{eq:shifted-saturated-moment}
 \sum_{n\in\Ssmooth(y)}\frac{(\log n)^\ell}{n^\sigma}
 =(1+o(1))\mu_y^\ell Z_y(\sigma).
\end{equation}
For $\ell=0$, the Euler-product identity holds without error.

To control the tail, shift $\sigma$ to $\sigma-1/\log_2X$.
The new parameters are $A+1$ and $B-1$.
The preceding moment calculation remains uniform because $B\to\infty$.
The prime number theorem, or partial summation of its elementary upper
bound, gives
\begin{align*}
 \log\frac{Z_y(\sigma-1/\log_2X)}{Z_y(\sigma)}
 &\leq\log Z_y(\sigma-1/\log_2X)\\
 &\ll\log_3X+\frac{\e^{A+1}}{A+1}
 =o\left(\frac{\log X}{\log_2X}\right).
\end{align*}
The first-order prime sum is bounded by the displayed
expression.  The higher-order terms are
$O(\sum_{p\leq y}p^{-1})=O(\log_3X)$.
For the last equality, $A+B=\tfrac12\log_2X$ implies
$\e^{A+1}\leq e\sqrt{\log X}$, while $A\to\infty$.
Equation~\eqref{eq:shifted-saturated-mean} also shows that the ratio of
the shifted mean to $\mu_y$ is bounded, since it is
$(1+o(1))eA/(A+1)$.
Applying \eqref{eq:shifted-saturated-moment} at the two arguments gives
\[
 \log
 \frac{\displaystyle\sum_{n\in\Ssmooth(y)}
        (\log n)^\ell n^{-\sigma+1/\log_2X}}
      {\displaystyle\sum_{n\in\Ssmooth(y)}
        (\log n)^\ell n^{-\sigma}}
 =o\left(\frac{\log X}{\log_2X}\right).
\]
For $n>X^\eta$ one has
$n^{-1/\log_2X}\leq X^{-\eta/\log_2X}$.
Using this inequality term by term and summing gives the relative
tail estimate
\begin{align*}
 &\frac{\displaystyle\sum_{\substack{n>X^\eta\\n\in\Ssmooth(y)}}
       (\log n)^\ell n^{-\sigma}}
      {\displaystyle\sum_{n\in\Ssmooth(y)}
       (\log n)^\ell n^{-\sigma}}\\
 &\quad\leq X^{-\eta/\log_2X}
 \frac{\displaystyle\sum_{n\in\Ssmooth(y)}
       (\log n)^\ell n^{-\sigma+1/\log_2X}}
      {\displaystyle\sum_{n\in\Ssmooth(y)}
       (\log n)^\ell n^{-\sigma}}\\
 &\quad=\exp\left\{-\frac{\eta\log X}{\log_2X}
       +o\left(\frac{\log X}{\log_2X}\right)\right\}=o(1).
\end{align*}
The sum in \eqref{eq:shifted-saturated-main} is
$1-o(1)$ times the full smooth sum.
The full sum tends to infinity uniformly.
The first-order prime sum gives $\log Z_y(\sigma)\gg\e^A/A$,
which tends to infinity because $A\geq\tfrac12\log_3X$.
The finite negative-square contribution and all the absolute errors
are therefore negligible relative to the main term.
Division by \eqref{eq:shifted-saturated-den}, followed by the comparison
of a nonnegative weighted average with the maximum, proves
\eqref{eq:shifted-saturated}.
\end{proof}

\begin{proof}[Proof of Theorem~\ref{thm5}]
When $1\leq A\leq\tfrac12\log_3X$, Proposition~\ref{prop:saturated}
gives the smooth sum with relative factor $1-o(1)$.
The growing-parameter assertion in Lemma~\ref{lem:friable} evaluates
that sum as
\[
 (1+o(1))D_\ell(A)(\log_2X)^{\ell+1},
\]
uniformly in this range.  Multiplication of the two relative estimates
proves \eqref{eq:transition-strong}.

Suppose now that $A>\tfrac12\log_3X$ and $B\to\infty$.
Proposition~\ref{prop:shifted-saturated} and
\eqref{eq:shifted-saturated-moment} give the lower bound
$(1-o(1))\mu_y^\ell Z_y(\sigma)$, where $y=\log X$ and
\[
 \mu_y=(1+o(1))\frac{\e^A\log_2X}{A}.
\]
Since $\log(1-x)^{-1}\geq x$ for $0<x<1$,
Lemma~\ref{lem:prime-sum} gives
\[
 \log Z_y(\sigma)
 \geq\sum_{p\leq y}p^{-\sigma}
 \geq\sigma\log_3X+\frac{\e^A}{A}+O(1).
\]
The hypothesis of that lemma holds because
$(1-\sigma)\log y=A\to\infty$.
Exponentiation yields
$Z_y(\sigma)\gg(\log_2X)^\sigma\exp(\e^A/A)$.
To compare the remaining factors with the stated bound, compute
\begin{align*}
 \log\frac{\mu_y^\ell(\log_2X)^\sigma}
 {\e^{-A}(\log_2X)^{\ell+1}}
 &=(\ell+1)A-\ell\log A
   -\frac{A\log_3X}{\log_2X}+O_\ell(1)\\
 &\longrightarrow+\infty.
\end{align*}
The conclusion is uniform because $\log_3X/\log_2X=o(1)$,
$\log A=o(A)$, and $A\geq\tfrac12\log_3X\to\infty$.
The product $\mu_y^\ell Z_y(\sigma)$ is at least a
constant depending only on $\ell$ times
\[
 \exp\left(\frac{\e^A}{A}-A\right)
 (\log_2X)^{\ell+1}.
\]
This proves \eqref{eq:transition-standard} and completes the proof.
\end{proof}

\end{document}